\documentclass[12pt,a4paper]{article}
\usepackage{}
\usepackage{mathrsfs}
\usepackage{amssymb}
\usepackage{amsmath}
\usepackage{amsfonts}
\usepackage{amstext}
\usepackage{amsthm}
\usepackage{graphicx}
\usepackage[top=0.9in, bottom=0.9in,left=0.8in,right=0.8in]{geometry}

\def \qed {\hfill \vrule height6pt width 6pt depth 0pt}

\begin{document}

\title{Bifurcation and one-sign solutions of the $p$-Laplacian involving a nonlinearity with zeros
\thanks{Research supported by NNSF of China (No. 11261052, 11401477), the Fundamental Research Funds for the Central Universities (No. DUT15RC(3)018)
and Scientific Research Project of the Higher Education Institutions of Gansu Province (No. 2014A-009).}}
\author{{\small Guowei Dai\thanks{Corresponding author. 
\newline
\text{\quad\,\,\, E-mail address}: daiguowei@dlut.edu.cn.%
}
} \\
{\small School of Mathematical Sciences, Dalian University of Technology, Dalian, 116024, PR China}\\
}
\date{}
\maketitle

\begin{abstract}
In this paper, we use bifurcation method to investigate the existence and multiplicity of one-sign solutions of the $p$-Laplacian involving a linear/superlinear nonlinearity with zeros. To do this, we first establish a bifurcation theorem from infinity for nonlinear operator equation with homogeneous operator. To deal with the superlinear case, we establish several topological results involving superior limit.
\end{abstract}\textbf{\ }

\numberwithin{equation}{section}

\numberwithin{equation}{section}

\section{Introduction}
\bigskip
\quad\, Consider the following problem
\begin{equation}\label{zeros}
\left\{
\begin{array}{ll}
-\Delta_p u=\lambda f(x,u)\,\, &\text{in}\,\, \Omega,\\
u=0  &\text{on}\,\,\partial\Omega,
\end{array}
\right.
\end{equation}
where $p\in(1,+\infty)$, $\Delta_pu = \text{div}\left(\vert \nabla u\vert^{p-2}\nabla u\right)$ is the $p$-Laplacian, $\lambda> 0$ is a real parameter, $\Omega$ is a bounded domain of $\mathbb{R}^N$ with smooth boundary $\partial\Omega$, and where $f$ is a nonlinearity with two
zeros which may vary in the $x$ variable. We will investigate the existence and multiplicity of one-sign solutions of problem (\ref{zeros}) involving a linear/superlinear growth nonlinearity by using of bifurcation and topological methods.

Problems with linear/superlinear nonlinearities at infinity have been extensively studied. For the Laplacian, see for example [\ref{ABC}, \ref{AH}, \ref{157}, \ref{CDG}, \ref{24}, \ref{FGU}, \ref{FLN}, \ref{25}, \ref{27}, \ref{GT}, \ref{22}, \ref{30}, \ref{LS}, \ref{Lions}, \ref{MD}, \ref{P}, \ref{29}]. For the $p$-Laplacian, see for
example [\ref{AAP}, \ref{DM}, \ref{PS}]. In most of these works, the nonlinearity is superlinear growth at infinity and satisfying the signum condition $f(x,s)s>0$ for $s\neq 0$.
The problem whose nonlinearity has nontrivial zeros has attracted a great deal of attention in history.
In [\ref{Lions}], this type of problems is considered for the Laplacian operator and a nonlinearity $f$ that is independent
on $x$, satisfying $f(0)\geq0$, $f(\beta) = 0$, and which is positive and superlinear for $s >\beta> 0$. Using topological degree arguments and under additional technical conditions which ensure a priori bounds, it is
shown that there exist two positive solutions of problem (\ref{zeros}). It is further shown that one solution
lies strictly below $\beta$, while the other has a maximum greater than $\beta$. This type of problems was also
studied in [\ref{Liu}], where again the existence of two positive solutions of problem (\ref{zeros}) was shown which improves the corresponding results of [\ref{Lions}] in some sense.
The characteristics
of the problem are quite different when the nonlinearity is linear growth at infinity and has nontrivial zeros. In [\ref{AH}], by using of Rabinowitz's bifurcation Theorem, Ambrosetti and Hess studied the global behavior of the
component of positive solutions of linear elliptic equation involving an asymptotically
nonlinearity with positive zero.
In [\ref{AAP}], Ambrosetti et al. studied the existence of branch of positive solutions for
asymptotically equidiffusive problem, which extends the corresponding ones of [\ref{AH}].
Delgado and Su\'{a}rez [\ref{DS}] studied the existence of positive solutions for some indefinite nonlinear eigenvalue problems involving
nonlinearities with positive zero.
In [\ref{M}], Ma studied the global
behavior of the components of nodal solutions of one-dimensional asymptotically linear eigenvalue problems
involving nonlinearities with two nontrivial zeros. Dai et al. [\ref{DML}] studied nodal solutions of quasilinear problems without signum condition.
In the above papers, the zeros of nonlinearity are constants. In [\ref{IMSU}], mainly by variational method,
the authors obtained some results involving the existence, multiplicity,
and the behavior with respect to $\lambda$ of positive solutions of problem (\ref{zeros}) under some suitable assumptions on $f$.

Motivated by the above papers, we study the existence and multiplicity of one-sign solutions of problem (\ref{zeros}) by bifurcation and topological methods. However, the bifurcation results of [\ref{D1}, \ref{D2}, \ref{L1}, \ref{R2}, \ref{R3}, \ref{R4}] cannot be applied directly to quasilinear problems, which are very important tools in dealing with semilinear boundary value problems. So we first study the bifurcation phenomenon from infinity for the following abstract operator equation
\begin{equation}\label{aoe}
u=L(\lambda) u+K(\lambda,u),\,\, u\in E,
\end{equation}
where $\lambda$ varies in $\mathbb{R}$, $E$ is a real Banach space
with norm $\Vert\cdot\Vert$, the map $\lambda\mapsto L(\lambda)$ is continuous. Moreover, we assume that $L(\cdot):E\rightarrow E$ is a homogeneous completely continuous operator and $K:\mathbb{R}\times E\rightarrow E$ is continuous with $K(\lambda,u)=o\left(\Vert u\Vert\right)$ near $u=\infty$ uniformly on bounded $\lambda$ intervals and $\Vert u\Vert^2 K\left(\lambda,u/\Vert u\Vert^2\right)$ is compact.
If $L(\lambda)=\lambda L$ and $L$ is linear, Rabinowitz [\ref{R3}] has shown that
all characteristic values of $L$ which are of odd multiplicity are bifurcation points from infinity. More precise, if $\mu$ is of odd multiplicity and
$\mathscr{S}$ denotes the set of solutions of equation (\ref{aoe}),
then $\mathscr{S}$ possesses an unbounded component which meets $(\mu,\infty)$ and satisfies some alternatives.
His method is to reduce the problem of bifurcation at infinity to bifurcation from the trivial solution axis which can be founded in [\ref{K}, \ref{R2}].
We also refer to [\ref{AM}, \ref{Bg}, \ref{K0}, \ref{LSW}, \ref{LSW1}, \ref{L1}, \ref{LM2}, \ref{Shi}, \ref{SW}] and their references for the bifurcation from the trivial solution. Amann [\ref{Amann}] established some bifurcation results from infinity in ordered Banach space by using this method.

As usually, $\mu$ is called an eigenvalue value of
\begin{equation}\label{qle}
u=L(\lambda) u,\,\, u\in E
\end{equation}
if there exists $0\not\equiv \varphi\in E$ such that
$\varphi={L}(\mu) \varphi$. Let $\Sigma$ denote the set of real eigenvalues of equation (\ref{qle}). If $\mu\in \Sigma$ has geometric multiplicity 1, we let
$E_0$ be a closed subspace of $E$ such that
\begin{equation}
E=\text{span}\left\{\varphi_\mu\right\}\oplus E_0,\nonumber
\end{equation}
where $\varphi_\mu$ is an eigenfunction corresponding to $\mu$ with $\left\Vert \varphi_\mu\right\Vert=1$. Let $B_r(0)=\left\{u\in E:\left\Vert u\right\Vert<r\right\}$. The first main result of this paper is the following theorem.
\\ \\
\textbf{Theorem 1.1.} \emph{If $\mu\in \Sigma$ is isolated, has geometric multiplicity 1 and such that}
\begin{equation}\label{jump}
\deg \left(I-L(\mu-\varepsilon), B_r(0),0 \right)\neq\deg \left(I-L(\mu+\varepsilon), B_r(0),0 \right)
\end{equation}
\emph{for any $\varepsilon>0$ small enough, then $\mathscr{\mathscr{S}}$ possesses two unbounded components $\mathscr{D}_\mu^+$ and $\mathscr{D}_\mu^-$  which meet $(\mu,\infty)$. Moreover if $\Lambda\subset \mathbb{R}$ is an interval such that $\Lambda\cap \Sigma =\{\mu\}$ and $\mathscr{M}$ is a neighborhood of $(\mu,\infty)$ whose projection on $\mathbb{R}$
lies in $\Lambda$ and whose projection on $E$ is bounded away from 0, then at least one of the
following three properties is satisfied by $\mathscr{D}_\mu^\nu$ for $\nu=+$ and $-$:}
\\

1$^o$. \emph{$\mathscr{D}_\mu^\nu-\mathscr{M}$ is bounded in $\mathbb{R}\times E$ in which case $\mathscr{D}_\mu^\nu-\mathscr{M}$
meets $\mathscr{R}=\{(\lambda,0):\lambda\in \mathbb{R}\}$,}
\\

2$^o$. \emph{$\mathscr{D}_\mu^\nu-\mathscr{M}$ is unbounded,}\\

3$^o$. \emph{contains a point $(\lambda,v)\in \mathbb{R}\times\left(E_0\setminus\{0\}\right)$.} \\

\emph{If 2$^o$ occurs and $\mathscr{D}_\mu^\nu-\mathscr{M}$ has a bounded projection on $\mathbb{R}$, then
$\mathscr{D}_\mu^\nu-\mathscr{M}$ meets $\left(\widehat{\mu},\infty\right)$ where $\mu\neq\widehat{\mu}\in \Sigma$. Moreover, there exists a neighborhood $\mathscr{N}\subset\mathscr{M}$ of $(\mu,\infty)$ such that
$(\lambda,u)\in \left(\mathscr{D}_\mu^\nu \cap \mathscr{N}\right)$ and $(\lambda,u)\neq (\mu,\infty)$
implies $(\lambda,u)=\left(\mu+o(1),\alpha \varphi_\mu+z\right)$ where $\alpha>0$ ($\alpha<0$) and $z=o( \alpha)$ at $\alpha=\infty$.}\\

In order to deal with the case of superlinear growth at infinity, we establish a Whyburn type limit theorem as follows.
\\ \\
\textbf{Theorem 1.2.} \emph{Let $X$ be a normal space and let $\left\{C_n\right\}$ be a sequence of unbounded connected subsets of $X$. Assume that:}
\\

(i) \emph{there exists $z^*\in \liminf_{n\rightarrow +\infty} C_n$ with $\left\Vert z^*\right\Vert=+\infty$;}

(ii) \emph{there exists a homeomorphism $T:X\rightarrow X$ such that $\left\Vert T\left(z^*\right)\right\Vert<+\infty$ and $\left\{T\left(C_n\right)\right\}$ be a sequence of unbounded connected subsets of $X$};

(iii) \emph{for every $R>0$, $\left(\cup_{n=1}^{+\infty} T \left(C_n\right)\right)\cap \overline{B}_R$ is a relatively compact set of $X$.}
\\

\noindent \emph{Then $D:=\limsup_{n\rightarrow +\infty}C_n$ is unbounded closed connected.}\\

Now, we are in the position to state the following five hypotheses on the nonlinearity $f$.\\

(H1) The function $f:\overline{\Omega}\times \mathbb{R}\rightarrow\mathbb{R}$ is continuous.

(H2) There exist two nontrivial functions $a$, $b\in W^{1,p}(\Omega)\cap  C^1\left(\overline{\Omega}\right)$ such that $b\leq 0\leq a$, $f(x,b)=f(x,a)=0$, $f(x,s)s>0$ for $x\in\Omega$ and $s\neq b$, $0$, $a$ which are not the solutions of problem (\ref{zeros}).

(H3) There exists a function $\alpha\in L^\infty(\Omega)$ with $\alpha\geq0$, $\not\equiv 0$ such that
\begin{equation}
\lim_{\vert s\vert\rightarrow 0}\frac{f(x,s)}{\varphi_p(s)}=\alpha(x)\,\,\text{uniformly in}\,\, x\in \Omega,\nonumber
\end{equation}
where $\varphi_p(s)=\vert s\vert^{p-2}s$.\\

(H4) There exists a function $\beta\in L^\infty(\Omega)$ with $\beta\geq0$, $\not\equiv 0$ such that
\begin{equation}
\lim_{\vert s\vert\rightarrow +\infty}\frac{f(x,s)}{\varphi_p(s)}=\beta(x)\,\,\text{uniformly in}\,\, x\in \Omega.\nonumber
\end{equation}

(H5) There is a $q\in\left(p,p^*\right)$ such that
\begin{equation}
\lim_{ \vert s\vert\rightarrow+\infty}\frac{f(x,s)}{\varphi_q(s)}=\beta(x)\nonumber
\end{equation}
uniformly in $x\in\Omega$, where $\beta$ is the same as (H4), and
\begin{equation}\label{pb}
p^*=\left\{
\begin{array}{ll}
\frac{Np}{N-p}\,\,&\text{if}\,\, p<N,\\
+\infty\,\,&\text{if}\,\,p\geq N.
\end{array}
\right.\nonumber
\end{equation}

Let $\lambda_{1,\alpha}$ denote the first eigenvalue for the following nonlinear eigenvalue problem
\begin{equation}\label{ep0}
\left\{
\begin{array}{ll}
-\text{div}\left(\vert\nabla u\vert^{p-2}\nabla u\right)=\lambda \alpha(x)\varphi_p(u)\,\,&\text{in}\,\,\Omega,\\
u=0&\text{on}\,\,\partial\Omega.
\end{array}
\right.
\end{equation}
It is well known that $\lambda_{1,\alpha}$ is simple, isolated and the associated eigenfunction $\varphi_{1,\alpha}$ with $\left\Vert \varphi_{1,\alpha}\right\Vert_{C^1\left(\overline{\Omega}\right)}=1$ has one sign in $\Omega$ (see for example [\ref{A}, \ref{DPM}, \ref{TTU}]).
The following theorem is the first result on the existence of one-sign solutions of problem (\ref{zeros}).
\\ \\
\textbf{Theorem 1.3.} \emph{Let (H1)--(H3) hold. Then problem (\ref{zeros}) has at least two solutions $u_{0}^+$ and
$u_{0}^-$ for every $\lambda\in\left(\lambda_{1,\alpha},+\infty\right)$, such that $0<u_{0}^+\leq a$ and $b\leq u_{0}^-<0$ in $\Omega$.}\\
\\
\textbf{Remark 1.1.} We do not require $f$ satisfying the monotonicity hypothesis ($M_1$) of [\ref{IMSU}] which is essential to use sub and supersolution techniques. Moreover, we also do not require $a>0$ being weakly $p$-superharmonic function.
Due to these reasons, the conclusions of Theorem 1.3 improves the corresponding ones of [\ref{IMSU}, Theorem 1.1] in some sense even when $f$ is nonnegative.\\
\\
\textbf{Remark 1.2.} In fact, $\left(\lambda_{1,\alpha},0\right)$ is the only one bifurcation point on $\mathbb{R}\times\{0\}$ of one-sign solutions of problem (\ref{zeros}) (see Section 3). So under hypotheses (H1)--(H3) any one-sign solution of problem (\ref{zeros}) lies in one component. That is to say, we find the range of all one-sign solutions. So it only needs to study the structure and formulation of the component to find one-sign solutions of problem (\ref{zeros}).
\\

Furthermore, we can get the exact multiplicity of one-sign solutions for problem (\ref{zeros}) under more strict hypotheses than those of Theorem 1.3.
\\ \\
\textbf{Theorem 1.4.} \emph{Let (H1)--(H3) and the following hypothesis hold,}\\

(H6) \emph{$f(x,s)\equiv 0$ for any $x\in\Omega$ and any $s\not\in[b(x),a(x)]$, $f(\cdot,s)$ is $C^1$ with respect to $s\in[b(\cdot),a(\cdot)]$ and such that $f(\cdot,s)/\varphi_p(s)$ is decreasing in $[0,a(\cdot)]$ and is increasing in $[b(\cdot),0]$.}\\

\noindent \emph{Then, } \\

(1) \emph{problem (\ref{zeros}) has exactly two solutions $u^+(\lambda,\cdot)$ and $u^-(\lambda,\cdot)$ for $\lambda\in \left(\lambda_{1,\alpha},+\infty\right)$, such that $0<u^+(\lambda,\cdot)\leq a(\cdot)$ and $b(\cdot)\leq u^-(\lambda,\cdot)<0$ in $\Omega$;}\\

(2) \emph{all one-sign solutions of problem (\ref{zeros}) lie on two smooth curves
\begin{equation}
\mathscr{C}^\pm=\left\{\left(\lambda,u^\pm(\lambda,\cdot)\right)\,:\,
\lambda\in\left(\lambda_{1,\alpha}, +\infty\right)\right\},\nonumber
\end{equation}
$\mathscr{C}^+$ and $\mathscr{C}^-$ join at $\left(\lambda_{1,\alpha},0\right)$;}\\

(3) \emph{$u^+(\lambda,\cdot)$ ($u^-(\lambda,\cdot)$) is increasing (decreasing) with respect to $\lambda$.}\\

Let $\lambda_{1,\beta}$ denote the first eigenvalue of problem (\ref{ep0}) with the weighted function $\beta$. By using of Theorem 1.1, the results involving the linear growth at infinity can be obtained as follows.
\\ \\
\textbf{Theorem 1.5.} \emph{Let (H1)--(H4) hold. Assume that $\lambda_{1,\alpha}>\lambda_{1,\beta}$. Then}\\

(i) \emph{if $\lambda\in\left(\lambda_{1,\beta},\lambda_{1,\alpha}\right]$, then problem (\ref{zeros}) has at least two solutions
$u_{\infty}^+$ and $u_{\infty}^-$, such that $u_{\infty}^+$ is positive in $\Omega$ and $u_{\infty}^-$ is negative in $\Omega$;}\\

(ii) \emph{if $\lambda\in\left(\lambda_{1,\alpha},+\infty\right)$, then problem (\ref{zeros}) has at least four solutions $u_{\infty}^+$,
$u_{\infty}^-$, $u_{0}^+$ and $u_{0}^-$ such that $u_{\infty}^+$ and $u_{0}^+$ ($\leq a$) are positive in $\Omega$; $u_{\infty}^-$ and $u_{0}^-$ ($\geq b$) are negative in $\Omega$.}
\\ \\
\textbf{Theorem 1.6.} \emph{Let (H1)--(H4) hold. Assume that $\lambda_{1,\alpha}<\lambda_{1,\beta}$. Then}\\

(i) \emph{if $\lambda\in\left(\lambda_{1,\alpha},\lambda_{1,\beta}\right]$, then problem (\ref{zeros}) has at least two solutions $u_{0}^+$ and
$u_{0}^-$, such that $u_{0}^+$ ($\leq a$) is positive in $\Omega$ and $u_{0}^-$ ($\geq b$)
is negative in $\Omega$;}\\

(ii) \emph{if $\lambda\in\left(\lambda_{1,\beta},+\infty\right)$, then problem (\ref{zeros}) has at least four solutions $u_{\infty}^+$,
$u_{\infty}^-$,
$u_{0}^+$ and $u_{0}^-$ such that $u_{\infty}^+$, $u_{0}^+$ ($\leq a$) are positive in $\Omega$; $u_{\infty}^-$, $u_{0}^-$ ($\geq b$) are negative in $\Omega$.}
\\
\\
\textbf{Remark 1.3.} Clearly, if $\lambda_{1,\beta}=\lambda_{1,\alpha}$ and $\lambda\in\left(\lambda_{1,\alpha},+\infty\right)$, then problem (\ref{zeros}) has at least four solutions $u_{\infty}^+$,
$u_{\infty}^-$,
$u_{0}^+$ and $u_{0}^-$ such that $u_{\infty}^+$, $u_{0}^+$ ($\leq a$) are positive in $\Omega$; $u_{\infty}^-$, $u_{0}^-$ ($\geq b$) are negative in $\Omega$ (see Section 4).
\\
\\
\textbf{Remark 1.4.} Due to the dependence on $x$ of $f$ and zeros, and since we are not requiring continuously differentiable property of $f$, these results of Theorem 1.5 and 1.6 improve, in a certain sense, those of [\ref{AH}, \ref{AAP}] even when $p = 2$.
\\

Our result involving the superlinear growth at infinity is the following theorem.
\\ \\
\textbf{Theorem 1.7.} \emph{Let (H1)--(H3) and (H5) hold. Then}\\

(i) \emph{if $\lambda\in\left(0,\lambda_{1,\alpha}\right]$, then problem (\ref{zeros}) has at least two solutions
$u_{\infty}^+$ and $u_{\infty}^-$, such that $u_{\infty}^+$ is positive in $\Omega$ and $u_{\infty}^-$ is negative in $\Omega$;}\\

(ii) \emph{if $\lambda\in\left(\lambda_{1,\alpha},+\infty\right)$, then problem (\ref{zeros}) has at least four solutions $u_{\infty}^+$,
$u_{\infty}^-$, $u_{0}^+$ and $u_{0}^-$ such that $u_{\infty}^+$, $u_{0}^+$ ($\leq a$) are positive in $\Omega$; $u_{\infty}^-$, $u_{0}^-$ ($\geq b$) are negative in $\Omega$.}\\
\\
\textbf{Remark 1.5.} We would like to point out that the hypotheses of Theorem 1.7 are weaker than those of [\ref{IMSU}, Theorem 1.2--1.4].
We do not require $f$ satisfying ($H_5$), ($H_6$), conditions ($a$)--($d$), and the monotonicity hypothesis ($M_2$) of [\ref{IMSU}].
The conditions ($a$)--($d$) were used to prove that $u_0^+<a$ in $\Omega$ which is fundamental to be able to obtain $u_{\infty}^+$. The hypothesis ($H_5$) or ($H_6$) is required to prove that $u_{\infty}^+$ cannot stay below
$a$. And ($M_2$) is necessary in [\ref{IMSU}] to use sub and supersolution techniques. We do not need these conditions or hypotheses because we use bifurcation method, and do not need $u_0^+<a$ and the fact of $u_{\infty}^+$ cannot stay below
$a$. We only care $u_0^+\not\equiv u_{\infty}^+$. Of course, it can be guaranteed by $u_0^+<a$ and the fact of $u_{\infty}^+$ cannot stay below $a$ but
these are not necessary. \\

Finally, we state results about the asymptotic behavior of solutions obtained in Theorem 1.7.
\\ \\
\textbf{Theorem 1.8.} \emph{Under hypotheses (H1)--(H3), and the followings:}\\

($\widetilde{H}$5) \emph{hypothesis (H5) holds with $q\in \left(p, p_*\right)$, where $p_*$ denotes the Serrin's exponent given by $p_* =
(N-1)p/(N-p)$ ($\infty$) if $N > p$ ($N \leq p$)},

(H7) \emph{$a(x)>0$, $b(x)<0$ in $\Omega$, and there exists $\gamma > 0$ such that $f(x, s)\geq \gamma\left\vert s-a(x)\right\vert^{q-1}$ for $s\geq a(x)$, and $f(x, s)\leq - \gamma\left\vert s-b(x)\right\vert^{q-1}$ for $s\leq b(x)$,}\\

\noindent \emph{if there exists an $\epsilon>0$ such that $\epsilon \varphi_{1,\alpha}\leq u_{0}^+$, $u_{\infty}^+$, and $u_{0}^-$, $u_{\infty}^-\leq-\epsilon \varphi_{1,\alpha}$ which come from Theorem 1.7, then $\lim_{\lambda\rightarrow +\infty}u_{0}^+=\lim_{\lambda\rightarrow +\infty}u_{\infty}^+=a$ and $\lim_{\lambda\rightarrow +\infty}u_{0}^-=\lim_{\lambda\rightarrow +\infty}u_{\infty}^-=b$ pointwise in $\Omega$.}
\\ \\
\textbf{Remark 1.6.} The asymptotic behavior of solutions when $\lambda\rightarrow 0$ or $\lambda\rightarrow \lambda_{1,\alpha}$ is the natural result of bifurcation method. We also would like to point out that the solutions obtained in this paper are the sense of $C^1$ weak solutions. By standard elliptic regularity theory (see [\ref{F}, \ref{FZ}, \ref{To}]), we know that any weak solution of problem (\ref{zeros}) belongs to $C^{1,\delta}(\overline{\Omega})$ with some constant $\delta\in(0,1)$ under the hypotheses of (H1) and (H4).\\
\\
\noindent\textbf{Remark 1.7.} We require $a$, $b\in C^1\left(\overline{\Omega}\right)$ because we need to use $\mathbb{R}\times \{a\}$
and $\mathbb{R}\times \{b\}$ to separate the bifurcation branches from the trivial solution axis and infinity. The reason to choose
$C_0^1\left(\overline{\Omega}\right)$ as our state space rather than $W_0^{1,p}$ is that we need to compare the norm $\Vert \cdot\Vert_\infty$
and the norm of state space, which is easy to realize in $C_0^1\left(\overline{\Omega}\right)$.
\\

The rest of this paper is arranged as follows. In Section 2, we show or present some necessary preliminary results on bifurcation and topology, and give
the proofs of Theorem 1.1 and Theorem 1.2. Section 3 is devoted to prove Theorem 1.3 and 1.4.
In Section 4, we study the bifurcation phenomenon from infinity of problem (\ref{zeros}) with linear growth nonlinearity at infinity and give the proofs of Theorem 1.5 and 1.6. In the last Section, we study the bifurcation phenomenon from infinity of problem (\ref{zeros}) with superlinear growth nonlinearity at infinity and give the proofs of Theorem 1.7 and 1.8.

\section{Preliminaries on bifurcation and topology}
\bigskip
\quad\,  In this section, we first study the bifurcation phenomenon from infinity for equation (\ref{aoe}) and prove Theorem 1.1.
Then we give the proof of Theorem 1.2 and establish several topological results involving superior limit, which we believe could be useful and interesting by themselves.
\subsection{Global bifurcation results}
\bigskip
\quad\, Suppose that $H:\mathbb{R}\times E\rightarrow E$ is completely continuous
and $H(\lambda,u)=o(\Vert u\Vert)$ at $u=0$ uniformly on bounded $\lambda$ intervals.
Let $\mathcal{S}$ be the closure of the set of nontrivial solution pairs of the following equation
\begin{equation}\label{oe0}
u=L(\lambda) u+H(\lambda,u),\,\, u\in E.
\end{equation}

\indent Firstly, we can get the following Rabinowitz type global bifurcation results.
\\ \\
\textbf{Lemma 2.1.} \emph{If $\mu\in \Sigma$ is isolated and such that (\ref{jump}) holds, then $(\mu,0)$ is a bifurcation point for equation (\ref{oe0}). Moreover, $\mathcal{S}$
possesses a maximal continuum $\mathscr{C}_\mu$ such that $(\mu,0)\in \mathscr{C}_\mu$ and $\mathscr{C}_\mu$ either}\\

(i) \emph{meets infinity in $\mathbb{R}\times E$, or}

(ii) \emph{meets $\left(\widehat{\mu},0\right)$, where $\mu\neq {\widehat{\mu}}\in \Sigma$.}
\\ \\
\textbf{Proof.} Let us suppose by contradiction that $\left(\mu, 0\right)$
is not a bifurcation point. Then
there exist $\varepsilon > 0$, $\rho_0> 0$ such that for $\left\vert \lambda-\mu\right\vert\leq2\varepsilon$ and $0<\rho < \rho_0$ there is no
nontrivial solution of the equation
\begin{equation}
u-L(\lambda) u-H(\lambda,u)=0\nonumber
\end{equation}
with $\Vert u\Vert=\rho$. From the invariance of the degree under a compact
homotopy we obtain that
\begin{equation}\label{edc}
\text{deg}\left(I-\Psi(\lambda,\cdot),B_\rho(0),0\right)\equiv constant
\end{equation}
for $\lambda\in\left[\mu-2\varepsilon,\mu+2\varepsilon\right]$, where $\Psi(\lambda,u)=L(\lambda) u+H(\lambda,u)$.

By taking $\varepsilon$ smaller if necessary, we can assume that there is no real eigenvalue value
of (\ref{qle}) in $\left(\mu,\mu+2\varepsilon\right]$. Take $\lambda=\mu+\varepsilon$.
We claim that the equation
\begin{equation}\label{es}
u- L(\lambda)u-sH(\lambda,u)=0
\end{equation}
has no solution $u$ with $\Vert u\Vert=\rho$ for every $s\in[0, 1]$ and $\rho$ sufficiently small.
Suppose on the contrary, let $\left\{u_n\right\}$ be the solution of equation (\ref{es}) with $\left\Vert u_n\right\Vert\rightarrow 0$
as $n\rightarrow+\infty$.
Let $v_n=u_n/\left\Vert u_n\right\Vert$. It follows that $v_n$ is a solution of the following problem
\begin{equation}\label{evs}
v= L(\lambda)v+s\frac{H\left(\lambda,u_n\right)}{\left\Vert u_n\right\Vert}.
\end{equation}
By equation (\ref{evs}), we obtain that for some convenient subsequence
$v_n\rightarrow v_0$ as $n\rightarrow+\infty$. Now $v_0$ verifies the following equation
\begin{equation}
v= L(\lambda)v\nonumber
\end{equation}
and $\left\Vert v_0\right\Vert = 1$. This implies that $\lambda=\mu+\varepsilon$ is an eigenvalue value of (\ref{qle}). This is a contradiction.

From the invariance of the degree under
homotopies, we obtain that
\begin{equation}\label{edFk}
\deg\left(I-\Psi(\mu+\varepsilon,\cdot), B_r(0),0\right)=\deg\left(I-L(\mu+\varepsilon), B_r(0),0\right).
\end{equation}
Similarly, we can find that
\begin{equation}\label{edFk1}
\deg\left(I-\Psi((\mu-\varepsilon),\cdot), B_r(0),0\right)=\deg\left(I-L(\mu-\varepsilon), B_r(0),0\right).
\end{equation}
Relations (\ref{edc}), (\ref{edFk}) and (\ref{edFk1}) imply that
\begin{equation}\nonumber
\deg\left(I-L(\mu+\varepsilon), B_r(0),0\right)=\deg\left(I-L(\mu-\varepsilon), B_r(0),0\right).
\end{equation}
This is a contradiction. Hence, $\left(\mu, 0\right)$ is a bifurcation point.

By standard arguments in global bifurcation theory (see [\ref{R2}]), we can show the existence of
a global branch of solutions of (\ref{oe0}) emanating from $\left(\mu, 0\right)$.\qed
\\

From now on we assume that $\mu$ has geometric multiplicity 1. Let $\mathbb{E}=\mathbb{R}\times E$.
Given any $\iota\in \mathbb{R}$
and $0 < s < +\infty$, we consider an
open neighborhood of $(\iota, 0)$ in $\mathbb{E}$ defined by
\begin{equation}
\mathbb{B}_s\left(\iota,0\right)=\left\{(\lambda,u)\in \mathbb{E}:\Vert u\Vert+\vert\lambda-\iota\vert<s\right\}.\nonumber
\end{equation}
According to the Hahn-Banach theorem,
there exists a linear functional $l\in E^*$, here $E^*$ denotes the dual space of $E$, such that
\begin{equation}
l\left(\varphi_\mu\right)=1\,\, \text{and}\,\, E_0=\{u\in E:l(u)=0\}.\nonumber
\end{equation}
For any $0 < \eta < 1$, define
\begin{equation}
K_\eta=\left\{(\lambda,u)\in \mathbb{E}: \vert l(u)\vert>\eta\Vert u\Vert\right\}.\nonumber
\end{equation}
Obviously, $K_\eta$ is an open subset of $\mathbb{E}$ consisting of two disjoint components $K_\eta^+$ and
$K_\eta^-$, where
\begin{equation}
K_\eta^+=\left\{(\lambda,u)\in \mathbb{E}:  l(u)>\eta\Vert u\Vert\right\},\nonumber
\end{equation}
\begin{equation}
K_\eta^-=\left\{(\lambda,u)\in \mathbb{E}:  l(u)<-\eta\Vert u\Vert\right\}.\nonumber
\end{equation}
Clearly, both $K_\eta^+$ and $K_\eta^-$ are convex cones, $K_\eta^-=-K_\eta^+$, and $\nu t\varphi_\mu\in K_\eta^\nu$ for
every $t > 0$, where $\nu\in\{+,-\}$.\\

Applying the similar method to prove [\ref{L1}, Lemma
6.4.1] with obvious changes, we may obtain the following result.
\\ \\
\textbf{Lemma 2.2.} \emph{For every $\eta\in(0, 1)$ there exists a number $\delta_0>0$ such that
for each $0<\delta<\delta_0$,
\begin{equation}
\left(\left( \mathcal{S}\setminus\left\{\left(\mu,0\right)\right\}\right)\cap \overline{\mathbb{B}}_\delta\left(\mu,0\right)\right)\subseteq K_\eta.\nonumber
\end{equation}
Moreover, for each
\begin{equation}
(\lambda,u)\in\left( \mathcal{S}\setminus\left\{\left(\mu,0\right)\right\}\right)\cap \left(\overline{\mathbb{B}}_\delta\left(\mu,0\right)\right),\nonumber
\end{equation}
there are $s\in\mathbb{R}$ and unique $y\in E_0$ such that
\begin{equation}
u=s\varphi_\mu+y\,\, \text{and} \,\, \vert s\vert>\eta\Vert u\Vert.\nonumber
\end{equation}
Furthermore, for these solutions $(\lambda,u)$,
\begin{equation}
\lambda=\mu+o(1)\,\, \text{and}\,\, y=o(s)\nonumber
\end{equation}
as $s\rightarrow 0$.}\\ \\
\textbf{Lemma 2.3.} \emph{If (\ref{jump}) holds, then $\mathscr{C}_\mu$ possesses a subcontinuum in each of the cones}
\begin{equation}
K_\eta^+\cup\left\{\left(\mu,0\right)\right\}\,\, \text{and}\,\, K_\eta^-\cup\left\{\left(\mu,0\right)\right\}\nonumber
\end{equation}
\emph{each of which meets $(\mu,0)$ and $\partial\overline{\mathbb{B}}_\varrho\left(\mu,0\right)$ for all $\varrho>0$ that is sufficiently small.}
\\ \\
\textbf{Proof.} According to Lemmata 2.1 and 2.2, the result is true for at least one of the cones. Suppose it is not true, for example, for $K_\eta^-$.
Then, no continuum
\begin{equation}
\widetilde{\mathscr{C}}\subseteq \left(K_\eta^-\cup\left\{\left(\mu,0\right)\right\}\right)\nonumber
\end{equation}
exists such that
\begin{equation}
\widetilde{\mathscr{C}}\cap K_\eta^-\cap\partial{\mathbb{B}}_\varrho\left(\mu,0\right)\neq\emptyset\nonumber
\end{equation}
for each $\varrho>0$ that is sufficiently small.

Define
\begin{equation}\label{e1.1}
\widehat{H}(\lambda,u)=\left\{
\begin{array}{ll}
H(\lambda,u)\,\,\,\quad\quad\quad\quad\quad\quad\quad\,\,\, &\text{if\,\,}l(u)< -\eta\Vert u\Vert,\\
\frac{-l(u)}{\eta \Vert u\Vert}H\left(\lambda,-\eta\Vert u\Vert\varphi_\mu+y\right)\,\,\, &\text{if\,\,}-\eta\Vert u\Vert\leq l(u)\leq0,\\
-H(\lambda,-u)\,\,\, \quad\quad\quad\quad\quad\quad &\text{if\,\,}l(u)>0
\end{array}
\right.\nonumber
\end{equation}
and
\begin{equation}
{\Phi}(\lambda,u)=u- L (\lambda)u+\widehat{H}(\lambda,u).\nonumber
\end{equation}
The new mapping $\widehat{H}(\lambda,u)$ satisfies the same continuity properties as ${H}(\lambda,u)$
and is an odd function of $u$. According to Lemma 2.1, the equation ${\Phi}(\lambda,u)=0$ possesses a component
$\widehat{\mathscr{C}}_\mu$ of nontrivial solutions emanating from $(\mu,0)$. By Lemma 2.2, there exists a constant $\delta_0>0$ such that
for each $\varrho\in\left(0,\delta_0\right)$,
\begin{equation}
\left(\widehat{\mathscr{C}}_\mu\cap {\mathbb{B}}_\varrho\left(\mu,0\right)\right)\subseteq \left(K_\eta\cup\{(\mu,0)\}\right).\nonumber
\end{equation}
Due to the homotopy invariance of the degree, there exists $\varrho_1\in\left(0,\delta_0\right)$ such that for each $\varrho\in\left(0,\varrho_1\right)$, we
have
\begin{equation}
\widehat{\mathscr{C}}_\mu\cap \partial{\mathbb{B}}_\varrho\left(\mu,0\right)\cap K_\eta\neq \emptyset.\nonumber
\end{equation}
Since ${\Phi}(\lambda,u)$ is odd in $u$, we have that
\begin{equation}
\widehat{\mathscr{C}}_\mu\cap K_\eta^+=\left\{(\lambda,-u):(\lambda,u)\in \left(\widehat{\mathscr{C}}_\mu\cap K_\eta^-\right)\right\}.\nonumber
\end{equation}
It follows that for each $0<\varrho<\varrho_1$,
\begin{equation}
\widehat{\mathscr{C}}_\mu\cap K_\eta^-\cap \partial{\mathbb{B}}_\varrho\left(\mu,0\right)\neq \emptyset.\nonumber
\end{equation}
This is a contradiction.\qed
\\

Let $\delta > 0$ be the constant from Lemma 2.2 and $\nu\in\{+,-\}$. For $0 <\varepsilon \leq\delta$
we define $\mathscr{D}_{\mu,\varepsilon}^\nu$ to be the component of $\left\{\left(\mu, 0\right)\right\}\cup
\left( \mathcal{S}\cap \overline{\mathbb{B}}_\varepsilon \cap K_\eta^\nu\right)$ containing
$\left(\mu, 0\right)$, $\mathscr{C}_{\mu,\varepsilon}^\nu$ to be the component of $\overline{\mathscr{C}_\mu
\setminus \mathscr{D}_{\mu,\varepsilon}^{-\nu}}$ containing $\left(\mu, 0\right)$, and $\mathscr{C}_\mu^\nu$ to be the closure of
$\cup_{0<\varepsilon\leq \delta}\mathscr{C}_{\mu,\varepsilon}^\nu$. Obviously, $\mathscr{C}_\mu^\nu$ is connected.
In view of Lemma 2.2, the definition of $\mathscr{C}_\mu^\nu$ is independent
from the choice of $\eta$ and $\mathscr{C}_\mu=\mathscr{C}_\mu^+\cup\mathscr{C}_\mu^-$.
\\ \\
\textbf{Lemma 2.4.} \emph{If the assumptions of Lemma 2.1 hold and $\mu$ has geometric multiplicity 1, $\mathscr{C}_\mu$ can be decomposed into two sub-continua $\mathscr{C}_\mu^+$ and $\mathscr{C}_\mu^-$, such that each of them either satisfies the alternatives of Lemma 2.1 or contains a point $(\lambda,v)\in \mathbb{R}\times\left(E_0\setminus\{0\}\right)$.}
\\ \\
\textbf{Proof.} Suppose to the contrary that $\mathscr{C}_\mu^-$ does not satisfy any of the alternatives of the statement. By argument
similar to that of Theorem 6.4.3 of [\ref{L1}], we can show that there exists $\eta_0\in(0,1)$ such that for each $\eta\in\left(0,\eta_0\right)$,
\begin{equation}\label{32}
\mathscr{C}_\mu^-\subseteq \left(K_\eta^-\cup\left\{\left(\mu,0\right)\right\}\right)
\end{equation}
Now, consider the equation
\begin{equation}
{\Phi}(\lambda,u)=0,\nonumber
\end{equation}
where ${\Phi}$ is the same as the proof of Lemma 2.3.
Lemma 2.1 shows that a continuum $\widetilde{\mathscr{C}}_\mu$
bifurcates from $\left(\mu,0\right)$. According to (\ref{32}) and the fact of ${\Phi}$ is odd in $u$,
we find that
\begin{equation}
\widetilde{\mathscr{C}}_\mu=\mathscr{C}_\mu^-\cup \left\{(\lambda,-u):(\lambda,u)\in \mathscr{C}_\mu^-\right\}.\nonumber
\end{equation}
It follows that $\mathscr{C}_\mu^-$ is unbounded.
A symmetric argument
shows that $\mathscr{C}_\mu^+$ also satisfies the alternatives of the statement.\qed
\\

Now, using Lemma 2.1 and the method of [\ref{R3}], we can get the following theorem of bifurcation from infinity.
\\ \\
\textbf{Proposition 2.1.} \emph{If $\mu\in \Sigma$ is isolated and such that (\ref{jump}) holds, $\mathscr{\mathscr{S}}$ possesses an unbounded component $\mathscr{D}_{\mu}$ which meets $(\mu,\infty)$. Moreover if $\Lambda\subset \mathbb{R}$ is an interval such that $\Lambda\cap \Sigma =\{\mu\}$ and $\mathscr{M}$ is a neighborhood of $(\mu,\infty)$ whose projection on $\mathbb{R}$
lies in $\Lambda$ and whose projection on $E$ is bounded away from 0, then either}
\\

1$^o$. \emph{$\mathscr{D}_\mu-\mathscr{M}$ is bounded in $\mathbb{R}\times E$ in which case $\mathscr{D}_\mu-\mathscr{M}$
meets $\mathscr{R}=\{(\lambda,0):\lambda\in \mathbb{R}\}$ or}
\\

2$^o$. \emph{$\mathscr{D}_\mu-\mathscr{M}$ is unbounded.}\\

\emph{If 2$^o$ occurs and $\mathscr{D}_\mu-\mathscr{M}$ has a bounded projection on $\mathbb{R}$, then
$\mathscr{D}_\mu-\mathscr{M}$ meets $\left(\widehat{\mu},\infty\right)$ where $\mu\neq\widehat{\mu}\in \Sigma$.}\\
\\
\textbf{Proof.} If $(\lambda,u)\in\mathscr{S}$ with $\Vert u\Vert\neq 0$, dividing equation (\ref{aoe}) by $\Vert u\Vert^2$ and
setting $w=u/\Vert u\Vert^2$ yield
\begin{equation}\label{oe1}
w=L(\lambda)w+H(\lambda,w),
\end{equation}
where $H(\lambda,w)=\Vert w\Vert^2 K\left(\lambda,w/\Vert w\Vert^2\right)$ for $w\neq 0$.
Extend $H$ to $w=0$ by $H(\lambda,0)=0$.
It is obvious that $(\lambda,0)$ is always the solution of equation (\ref{oe1}).
It is not difficult to show that $H:\mathbb{R}\times E\rightarrow E$ is completely continuous
and $H(\lambda,w)=o\left(\Vert w\Vert\right)$ at $w=0$ uniformly on bounded $\lambda$ intervals.
Now applying Lemma 2.1 to equation (\ref{oe1}), we have the component $\mathscr{C}_\mu$ of
$\mathcal{S}$, containing $(\mu,0)$ and satisfying the alternatives of Lemma 2.1.
Under the inversion $w\rightarrow w/\Vert w\Vert^2=u$, $\mathscr{C}_\mu\rightarrow \mathscr{D}_\mu$
satisfying equation (\ref{aoe}). Clearly, $\mathscr{D}_\mu$ satisfies the conclusions of this theorem.\qed
\\

Furthermore, if $\mu$ has geometric multiplicity 1, we can use Lemma 2.4 instead of Lemma 2.1 in the proof of Proposition 2.1 to get Theorem 1.1.\\
\\
\textbf{Proof of Theorem 1.1.} Clearly, if $(\lambda,v)\in \mathscr{C}_\mu^\nu\cap \left(\mathbb{R}\times\left(E_0\setminus\{0\}\right)\right)$, then $(\lambda,u)\in \mathscr{D}_\mu^\nu$ with $u=v/\Vert v\Vert^2$. So we only need to show that there exists a neighborhood $\mathscr{N}\subset\mathscr{M}$ of $(\mu,\infty)$ such that
$(\lambda,u)\in \left(\mathscr{D}_\mu^\nu \cap \mathscr{N}\right)$ and $(\lambda,u)\neq (\mu,\infty)$
implies $(\lambda,u)=(\mu+o(1),\alpha \varphi_\mu+z)$.
Next, we only prove the case of $\nu=+$ since the proof of $\nu=-$ is similar. Clearly, the inversion $w\rightarrow w/\Vert w\Vert^2=u$ turns
$(\mu,0)$ into $(\mu,\infty)$.

By Lemma 2.2, there exists a bounded neighborhood $\mathscr{O}$ of $(\mu,0)$ such that for each
\begin{equation}
(\lambda,w)\in\left(\left( \mathscr{C}_\mu^+\setminus\left\{\left(\mu,0\right)\right\}\right)\cap \mathscr{O}\right),\nonumber
\end{equation}
there are (unique) $0\neq s\in\mathbb{R}^+$ ($\mathbb{R}^-$) and $y\in E_0$ such that $w=s\varphi_\mu+y$.
Furthermore, for these solutions $(\lambda,w)$, $\lambda=\mu+o(1)$ and $y=o(s)$ as $s\rightarrow 0$.
By the inversion $w\rightarrow w/\Vert w\Vert^2=u$, $\left( \mathscr{C}_\mu^+\setminus\left\{\left(\mu,0\right)\right\}\right)\cap \mathscr{O}$ is translated to a deleted neighborhood $\mathscr{N}^o$ of $(\mu,\infty)$. For these solutions $(\lambda,u)$, one has
\begin{equation}
u=\alpha \varphi_\mu+z,\nonumber
\end{equation}
where $\alpha=s/\Vert w\Vert^2$ and $z=y/\Vert w\Vert^2$.
It is sufficient to show $z=o(\alpha)$ at $\alpha=\infty$.
It is easy to verify that
\begin{equation}
\lim_{s\rightarrow 0} \frac{\Vert w\Vert}{\vert s\vert}=1.\nonumber
\end{equation}
It follows that $\lim_{s\rightarrow 0}  \alpha=\infty$ and $\lim_{\alpha\rightarrow \infty}  s=0$. Obviously, one has that
\begin{equation}
\lim_{s\rightarrow 0}\frac{z}{\alpha}=\lim_{s\rightarrow 0}\frac{y}{s}=0.\nonumber
\end{equation}
So we get that
\begin{equation}
\lim_{\alpha\rightarrow\infty}\frac{z}{\alpha}=0.\nonumber
\end{equation}
Therefore, the conclusions hold by taking $\mathscr{N}:=\left(\left(\mathscr{N}^o\cup \{(\mu,\infty)\}\right)\cap \mathscr{M}\right)$.\qed

\subsection{Topological results involving superior limit}
\bigskip
\quad\, The Whyburn's limit theorem [\ref{Whyburn}, Theorem 9.1] is an important tool in the study of differential equations theory, see for example, [\ref{ACD}, \ref{AT}, \ref{ADT}, \ref{Dai}] and references cited therein. However, if the collection of the infinite sequence of sets
is unbounded, the Whyburn's limit theorem cannot be used directly because the collection may not be relatively compact.

For the convenience of the reader, we recall the definitions of superior limit and inferior limit here. Let $G$ be any infinite
collection of point sets, not necessarily different. The set of all points
$x$ of our space $G$ such that every neighborhood of $x$ contains points of
infinitely many sets of $G$ is called the superior limit or limit superior
of $G$ and is written $\limsup G$. The set of all points $y$ such that every
neighborhood of $y$ contains points of all but a finite number of the sets
of $G$ is called the inferior limit or limit inferior of $G$ and is written $\liminf G$.\\

To prove Theorem 1.2, we need the following topological lemma.
\\ \\
\textbf{Lemma 2.5.} \emph{Let $\mathcal{X}$ be a normal space and let $\left\{\mathcal{C}_n\right\}$ be a sequence of unbounded connected subsets of $\mathcal{X}$. Assume that:}
\\

(i) \emph{there exists $z^*\in \liminf_{n\rightarrow +\infty} \mathcal{C}_n$ with $\left\Vert z^*\right\Vert<+\infty$;}

(ii) \emph{for every $R>0$, $\left(\cup_{n=1}^{+\infty} \mathcal{C}_n\right)\cap \overline{B}_R$ is a relatively compact set of $\mathcal{X}$}\\

\noindent \emph{Then $\mathcal{D}:=\limsup_{n\rightarrow +\infty}\mathcal{C}_n$ is unbounded closed connected.}
\\ \\
\textbf{Proof.} Let $X_R=\mathcal{X}\cap B_R$ for any $R>0$. Then $X_R$ is a metric subspace under the induced topology of $\mathcal{X}$. Let $A_n=\mathcal{C}_n\cap B_R$. Clearly, we have
$\cup_{n=1}^{+\infty} A_n=\left(\cup_{n=1}^{+\infty} \mathcal{C}_n\right)\cap B_R$.
So $\cup_{n=1}^{+\infty} A_n$ is relatively compact in $X_R$. Furthermore, $z^*\in \liminf_{n\rightarrow +\infty} \mathcal{C}_n$ implies that every
neighborhood $U\left(z^*\right)$ of $z^*$ contains points of all but a finite number of the sets
of $\left\{\mathcal{C}_n\right\}$. So there exists a positive integer $N$ such that for $n >N$, $U\left(z^*\right)\cap \mathcal{C}_n\neq\emptyset$.
Since $\left\Vert z^*\right\Vert<+\infty$, we can take $R>0$ large enough such that $U\left(z^*\right)\subseteq B_R$. Thus, $U\left(z^*\right)\cap A_n=U\left(z^*\right)\cap \mathcal{C}_n\neq\emptyset$ for $n >N$.
So we have $z^*\in \liminf_{n\rightarrow +\infty} A_n$.
By Theorem 9.1 of [\ref{Whyburn}], it follows that $A=\limsup_{n\rightarrow +\infty}A_n$ is connected in $X_R$.

We claim that $B:=\left(\limsup_{n\rightarrow +\infty}\mathcal{C}_n\right)\cap B_R=A$.
For $x\in A$, then any neighborhood $V$ in $X_R$ of $x$ contains contains points of
infinitely many sets of $\left\{\mathcal{C}_n\cap B_R\right\}$. So there exist $x_{n_i}\in \mathcal{C}_{n_i}\cap B_R$ such that
$x_{n_i}\rightarrow x$ as $i\rightarrow+\infty$. It follows that $x\in B_R$ and $x\in \limsup_{n\rightarrow +\infty}\mathcal{C}_n$, i.e., $x\in B$.
Conversely, if $x\in B$, any neighborhood $V$ in $X_R$ of $x$ contains a point $z$ of $\limsup_{n\rightarrow +\infty}\mathcal{C}_n$ and thus $V$, a neighborhood of $z$, contains points of infinitely many of the sets of $\mathcal{C}_n\cap B_R$.
It follows that $x\in A$. Hence, $B$ is connected. By the arbitrary of $R$, we get that $\mathcal{D}$ is connected.
From [\ref{Whyburn}], we know that $\mathcal{D}$ is closed.

Next, we show that $\mathcal{D}$ is unbounded.
Suppose on the contrary that $\mathcal{D}$ is bounded. It is easy to see that $\mathcal{D}$ is a compact set of $\mathcal{X}$ by (ii) and the fact of $z^*\in \mathcal{D}$.
Let $U_\delta$ be a $\delta$-neighborhood of $\mathcal{D}$. So we have that
\begin{equation}\label{Why}
\partial U_\delta\cap \mathcal{D}=\emptyset.
\end{equation}
By (i) and the connectedness of $\mathcal{C}_n$, there exists an integer $N_0 >0$, such that for all
$n>N_0$, $\mathcal{C}_n\cap \partial U_\delta\neq \emptyset$. Take $y_n \in \mathcal{C}_n\cap \partial U_\delta$, then $\left\{y_n:n>N_0\right\}$ is
a relatively compact subset of $\mathcal{X}$, so there exist $y^*\in \partial U_\delta$ and a subsequence $\left\{y_{n_k}\right\}$ such
that $y_{n_k}\rightarrow y^*$. The definition of superior limit shows that $y^*\in \mathcal{D}$. Therefore, $y^*\in \partial U_\delta\cap \mathcal{D}$. However, this
contradicts (\ref{Why}). Therefore, $\mathcal{D}$ is unbounded.\qed
\\ \\
\textbf{Proof of Theorem 1.2.} We first claim that $T\left(z^*\right)\in \liminf_{n\rightarrow +\infty} T\left(C_n\right)$.
By (i), we have that $T\left(z^*\right)\in T\left(\liminf_{n\rightarrow +\infty} C_n\right)$. It is enough to show that
\begin{equation}
T\left(\liminf_{n\rightarrow +\infty} C_n\right)=\liminf_{n\rightarrow +\infty} T\left(C_n\right).\nonumber
\end{equation}
For any $y\in \liminf_{n\rightarrow +\infty} T\left(C_n\right)$, let $x=T^{-1}(y)$. For any neighborhood $V$ in $X$ of $x$,
$T(V)$ is a neighborhood of $y=T(x)$. So $T(V)$ contains points of all but a finite number of the sets $\left\{T\left(C_n\right)\right\}$. Thus there exists $m\in \mathbb{N}$ such that $y_n\in T(V)$ for any $y_n\in T\left(C_n\right)$ and $n>m$.
Let $x_n=T^{-1}\left(y_n\right)$ for any $n>m$. Then obviously one has $x_n\in \left(C_n\cap V\right)$ for any $n>m$.
It follows that $x\in \liminf_{n\rightarrow +\infty} C_n$. Therefore, $y\in T\left(\liminf_{n\rightarrow +\infty} C_n\right)$.
Conversely, for any $y\in T\left(\liminf_{n\rightarrow +\infty} C_n\right)$, let $x\in \liminf_{n\rightarrow +\infty} C_n$ such that $y=T(x)$.
For any neighborhood $V$ in $X$ of $y$, $T^{-1}(V)$ is a neighborhood of $x=T^{-1}(y)$.
So $T^{-1}(V)$ contains points of all but a finite number of the sets $\left\{C_n\right\}$. Thus there exists $m\in \mathbb{N}$ such that $x_n\in T^{-1}(V)$ for any $x_n\in C_n$ and $n>m$. Let $y_n=T\left(x_n\right)$ for any $n>m$. Then obviously one has $y_n\in \left(T\left(C_n\right)\cap V\right)$ for any $n>m$.
It follows that $y\in \liminf_{n\rightarrow +\infty} T\left(C_n\right)$.

Now Lemma 2.5 shows that
$\widetilde{D}=\limsup_{n\rightarrow +\infty} T\left(C_n\right)$ is unbounded closed connected. Let $D=\limsup_{n\rightarrow +\infty}C_n$. Clearly, $z^*\in D$. The unboundedness of $z^*$ shows that $D$ is unbounded. From [\ref{Whyburn}], we know that $D$ is closed.
So we only need to show that $D$ is connected.

We next show that $T(D)=\widetilde{D}$.
For $y\in \widetilde{D}$, then any neighborhood $V$ in $X$ of $y$ contains points of
infinitely many sets of $\left\{D_n\right\}$. So there exist $y_{n_i}\in D_{n_i}$ such that
$y_{n_i}\rightarrow y$ as $i\rightarrow+\infty$. It follows that $y\in \limsup_{n\rightarrow +\infty}D_n$.
Let $x_{n_i}\in C_{n_i}$ such that $y_{n_i}=T\left(x_{n_i}\right)$ and $x=T^{-1}y$. Then one has that
$x_{n_i}=T^{-1}\left(y_{n_i}\right)\rightarrow T^{-1}y=x$. So $x\in D$, i.e., $y\in T\left(D\right)$.
Conversely, for any $y\in T(D)$, let $x\in D$ such that $y=T(x)$. Then any neighborhood $V$ in $X$ of $x$ contains contains points of
infinitely many sets of $\left\{C_n\right\}$. So there exist $x_{n_i}\in C_{n_i}$ such that
$x_{n_i}\rightarrow x$ as $i\rightarrow+\infty$. Let $y_{n_i}\in D_{n_i}$ such that $y_{n_i}=T\left(x_{n_i}\right)$.
It follows that $y_{n_i}\rightarrow y$. So we have $y\in \widetilde{D}$.

Finally, we show the connectedness of $D$. Suppose on the contrary that $D$ is not bounded connected.
So there exist two subsets $A$ and $B$ of $D$ such that $A\cap \overline{B}=B\cap \overline{A}=\emptyset$.
It is not difficult to verify that $T(A)\cap T\left(\overline{B}\right)=T(B)\cap T\left(\overline{A}\right)=\emptyset$.
That is to say $T(D)$ is not connected which contradicts the connectedness of $\widetilde{D}$.\qed\\

Note that the unboundedness of $\left\{T\left(C_n\right)\right\}$ is not necessary in the conclusion of Theorem 1.2. In fact,
by Theorem 9.1 of [\ref{Whyburn}], in view of the proof of Theorem 1.2, we can immediately get the following corollary.
\\ \\
\textbf{Corollary 2.1.} \emph{Let $X$ be a normal space and let $\left\{C_n\right\}$ be a sequence of unbounded connected subsets of $X$. Assume that:}
\\

(i) \emph{there exists $z^*\in \liminf_{n\rightarrow +\infty} C_n$ with $\left\Vert z^*\right\Vert=+\infty$;}

(ii) \emph{there exists a homeomorphism $T:X\rightarrow X$ such that $\left\Vert T\left(z^*\right)\right\Vert<+\infty$};

(iii) \emph{$\cup_{n=1}^{+\infty} T \left(C_n\right)$ is a relatively compact set of $X$.}\\

\noindent \emph{Then $D:=\limsup_{n\rightarrow +\infty}C_n$ is unbounded closed connected.}\\

The following lemma describes a relation of its superior limit and elements of a sequence of sets, which is very useful
to understand superior limit via element because superior limit may be more obscure than concrete element.
\\ \\
\textbf{Lemma 2.6} (see [\ref{Whyburn}]). \emph{Let $(X,\rho)$ be a metric space. If $\left\{A_i\right\}_{i\in \mathbb{N}}$ is a sequence of sets of $X$ whose limit superior
is $L$ and the sum of whose elements is a conditionally compact set, then for each $\epsilon> 0$ there exists an $m$ such that for every $n > m$, $A_n\subset V_\epsilon(L)$, where $V_\epsilon(L)$ denotes the set of all points $p$ with
$\rho(p, x) < \epsilon$ for any $x\in L$.}\\

In concrete problems, the collection of the infinite sequence of sets may not be relatively (conditionally) compact.
So the following two propositions are needed, which will play a key role in the analysis of the structure of bifurcation branch from infinity (see Section 5).
\\ \\
\textbf{Proposition 2.2.} \emph{If $\left\{C_i\right\}_{i\in \mathbb{N}}$ is a sequence of sets whose limit superior
is $L$ and for every $R>0$, $\left(\cup_{i=1}^{+\infty} C_i\right)\cap B_R$ is a relatively compact set, then for each $\epsilon> 0$ there exists an $m$ such that for every $n > m$, $C_n\subset V_\epsilon(L)$.}
\\
\\
\textbf{Proof.} Let $X_R=X\cap B_R$ for any $R>0$. Then $X_R$ is a metric subspace under the induced topology of $X$. Let $A_n=C_n\cap B_R$. Clearly, we have
$\cup_{n=1}^{+\infty} A_n=\left(\cup_{n=1}^{+\infty} C_n\right)\cap B_R$.
So $\cup_{n=1}^{+\infty} A_n$ is relatively compact in $X_R$.

Let $A=\limsup_{n\rightarrow +\infty}A_n$. We claim that $A=L\cap B_R:=B$.
For any $x\in A$, then any neighborhood $V$ in $X_R$ of $x$ contains points of
infinitely many sets of $\left\{C_n\cap B_R\right\}$. So there exist $x_{n_i}\in C_{n_i}\cap B_R$ such that
$x_{n_i}\rightarrow x$ as $i\rightarrow+\infty$. It follows that $x\in B_R$ and $x\in \limsup_{n\rightarrow +\infty}C_n$, i.e., $x\in B$.
Conversely, if $x\in B$, any neighborhood $V$ in $X_R$ of $x$ contains a point $z$ of $\limsup_{n\rightarrow +\infty}C_n$ and thus $V$, a neighborhood of $z$, contains points of infinitely many of the sets of $C_n\cap B_R$.
It follows that $x\in A$.

Now by Lemma 2.6, for each $\epsilon> 0$ there exists an $m$ such that for every $n > m$, $A_n\subset V_\epsilon(B)$ in $X_R$.
Clearly, one has $V_\epsilon(B)\subset V_\epsilon(L)$. Thus we have $A_n\subset V_\epsilon(L)$ in $X$.
By the arbitrary of $R$, we get that $C_n\subset V_\epsilon(L)$ in $X$ for every $n > m$.\qed
\\ \\
\textbf{Proposition 2.3.} \emph{If $\left\{C_i\right\}_{i\in \mathbb{N}}$ is a sequence of sets whose limit superior
is $L$ and there exists a homeomorphism $T:X\rightarrow X$ such that for every $R>0$, $\left(\cup_{i=1}^{+\infty} T\left(C_i\right)\right)\cap B_R$ is a relatively compact set, then for each $\epsilon> 0$ there exists an $m$ such that for every $n > m$, $C_n\subset V_\epsilon(L)$.}
\\
\\
\textbf{Proof.} By an argument similar to that of Theorem 1.2, we can obtain that $\limsup_{n\rightarrow +\infty}T\left(C_n\right)=T(L)$.
Applying Proposition 2.2 to $\left\{T\left(C_i\right)\right\}_{i\in \mathbb{N}}$, we get that for each $\widetilde{\epsilon}> 0$ there exists an $m$ such that for every $n > m$, $T\left(C_n\right)\subset V_{\widetilde{\epsilon}}(T(L))$. For any fixed $x_0\in C_n$, we have $T\left(x_0\right)\in T\left(C_n\right)$. Clearly, one has $T\left(x_0\right)\in V_{\widetilde{\epsilon}}(T(L))$.
It means that for any $z\in L$, $\left\Vert T\left(x_0\right)-T(z)\right\Vert<\widetilde{\epsilon}$.
It follows from the continuity of $T^{-1}$ that $\left\Vert x_0-z\right\Vert=\left\Vert T^{-1}\left(T\left(x_0\right)\right)-T^{-1}(T(z))\right\Vert<\epsilon$ for each $\epsilon> 0$. By the arbitrary of $x_0$, we get that $C_n\subset V_\epsilon(L)$ for every $n > m$.\qed

\section{Bifurcation from trivial solution axis}
\bigskip
\quad\, The aim of this section is to prove Theorem 1.3 and 1.4. Firstly, we define
\begin{equation}\label{e1.1}
\widetilde{f}(x,s)=\left\{
\begin{array}{ll}
f(x,s)\,\, &\text{if\,\,} b(x)\leq s\leq a(x),\\
0 &\text{otherwise}
\end{array}
\right.\nonumber
\end{equation}
and consider the following problem
\begin{equation}\label{ap}
\left\{
\begin{array}{ll}
-\Delta_p u=\lambda \widetilde{f}(x,u)\,\, &\text{in}\,\, \Omega,\\
u=0  &\text{on}\,\,\partial\Omega.
\end{array}
\right.
\end{equation}
Set
\begin{equation}
\mathbb{S}^+:=\left\{u\in C^{1,\delta}\left(\overline{\Omega}\right):u> 0\, \text{\,in\,\,}\Omega\right\}\nonumber
\end{equation}
and
\begin{equation}
\mathbb{S}^-:=\left\{u\in C^{1,\delta}\left(\overline{\Omega}\right):u< 0\, \text{\,in\,\,}\Omega\right\}.\nonumber
\end{equation}
Let $E=\left\{u\in C^1(\overline{\Omega}):u=0\,\,\text{on}\,\,\partial\Omega\right\}$ with the usual norm $\Vert u\Vert=\max_{\overline{\Omega}}\vert u\vert+\max_{\overline{\Omega}}\left\vert \nabla u\right\vert$.

Consider the following auxiliary problem
\begin{equation}\label{eh}
\left\{
\begin{array}{ll}
-\text{div}\left(\varphi_p(\nabla u)\right)=h\,\,  &\text{in}\ \ \Omega,\\
u=0&\text{on}\,\,\partial\Omega
\end{array}
\right.
\end{equation}
for a given $h\in L^{r/(r-1)}(\Omega)$, where $r\in\left(1,p^*\right)$.
By a solution of problem (\ref{eh}) we understand that $u\in W_0^{1,p}(\Omega)$ satisfying problem (\ref{eh}) in
the weak sense, i.e., such that
\begin{equation}\label{e1.1}
\int_\Omega \vert \nabla u\vert^{p-2}\nabla u\nabla v\,dx=\int_\Omega hv\,dx\nonumber
\end{equation}
for all $v\in W_0^{1,p}(\Omega)$.
We have known that for every given $h\in L^{r/(r-1)}(\Omega)$ there is a unique solution
$u$ to problem (\ref{eh}) (see [\ref{DPM}]). Let $R_p(h)$ denote the unique solution to problem (\ref{eh})
for a given $h\in L^{r/(r-1)}(\Omega)$. It is well known that $R_p:L^\infty(\Omega)\rightarrow E$ is completely continuous (see [\ref{F}, \ref{F1}, \ref{Li}]).
The definition of $\widetilde{f}$ implies that $u\in C^{1,\delta}\left(\overline{\Omega}\right)$ with some constant $\delta\in(0,1)$ for every weak solution $u$ of
problem (\ref{ap}).\\

Next we study the bifurcation phenomenon of problem (\ref{ap}) from the trivial solution axis.
\\ \\
\textbf{Proposition 3.1.} \emph{Under the assumptions of (H1)--(H3), the pair $\left(\lambda_{1,\alpha},0\right)$ is a bifurcation point of problem (\ref{ap}). Moreover, there exist two
unbounded continua of the set of nontrivial solutions of
problem (\ref{ap}) in $\mathbb{R}\times E$ bifurcating
from $\left(\lambda_{1,\alpha},0\right)$, $\mathscr{C}^+$ and $\mathscr{C}^-$, such that $\mathscr{C}^\nu\subseteq\left(\left\{\left(\lambda_{1,\alpha},0\right)\right\}\cup\left(\mathbb{R}\times \mathbb{S}^\nu\right)\right)$. }
\\ \\
\textbf{Proof.} Let $\xi:\Omega\times\mathbb{R}\rightarrow \mathbb{R}$ be such that
\begin{equation}
\widetilde{f}(x,s)=\alpha(x)\varphi_p(s)+\xi(x,s)\nonumber
\end{equation}
with
\begin{equation}
\lim_{\vert s\vert\rightarrow0}\frac{\xi(x,s)}{\varphi_p(s)}=0.\nonumber
\end{equation}
Let us consider
\begin{equation}\label{tp}
\left\{
\begin{array}{ll}
-\Delta_p u=\lambda \alpha(x)\varphi_p(u)+\lambda\xi(x,u)\,\, &\text{in}\,\, \Omega,\\
u=0  &\text{on}\,\,\partial\Omega
\end{array}
\right.
\end{equation}
as a bifurcation problem from the trivial solution axis.

Now problem (\ref{tp}) can be equivalently written as
\begin{equation}\label{qop}
u=R_p\left(\lambda K_pu+H(\lambda,u)\right),
\end{equation}
where $K_p(u)=\alpha(\cdot)\varphi_p(u)$, $H(\lambda,\cdot)$ denotes the usual Nemitsky operator associated with $\lambda\xi$.
Let
$L(\lambda)u=R_p\left(\lambda K_pu\right)$ and $\widehat{H}(\lambda,u)=R_p\left(\lambda K_pu+H(\lambda,u)\right)-L(\lambda)u$.
Then it is easy to see that $L(\cdot):E\rightarrow E$ is homogeneous completely continuous.
From conditions (H1)--(H3),
the right hand side of equation (\ref{qop}) defines a completely continuous operator from $\mathbb{R}\times E$ into $E$.
So $\widehat{H}:\mathbb{R}\times E\rightarrow E$ is completely continuous.

Let
\begin{equation}
\widetilde{\xi}(x,u)=\max_{0\leq \vert s\vert\leq u}\vert \xi(x,s)\vert\,\,
\text{for any}\,\, x\in \Omega\nonumber
\end{equation}
then $\widetilde{\xi}$ is nondecreasing with respect to $u$ and
\begin{equation}\label{h0}
\lim_{ u\rightarrow 0^+}\frac{\widetilde{\xi}(x,u)}{
u^{p-1}}=0.
\end{equation}
Further it follows from (\ref{h0}) that
\begin{equation}
\left\vert\frac{\xi(x,u)}{\Vert u\Vert^{p-1}}\right\vert \leq\frac{
\widetilde{\xi}(x,\vert u\vert)}{\Vert u\Vert^{p-1}} \leq \frac{
\widetilde{\xi}\left(x,\Vert u\Vert_\infty\right)}{\Vert u\Vert^{p-1}} \leq
\frac{ \widetilde{\xi}(x,\Vert u\Vert)}{\Vert
u\Vert^{p-1}}\rightarrow0\,\,  \text{as}\,\, \Vert
u\Vert\rightarrow 0\nonumber
\end{equation}
uniformly in $x\in\Omega$, where $\Vert u\Vert_\infty=\max_{x\in \overline{\Omega}} \vert u(x)\vert$. It follows that $\widehat{H}=o(\Vert u\Vert)$ near $u=0$ uniformly on bounded $\lambda$ intervals.
By an argument similar to that of [\ref{DPM}, Proposition 2.2], we can get that
\begin{equation}
\deg_E\left(I-L(\lambda), B_r(0),0 \right)=\left\{
\begin{array}{ll}
1 \,\, &\text{if}\,\, 0<\lambda<\lambda_{1,\alpha},\\
-1\,\,&\text{if}\,\,\lambda_{1,\alpha}<\lambda<\lambda_{2,\alpha},
\end{array}
\right.\nonumber
\end{equation}
where $\lambda_{2,\alpha}$ denotes the second eigenvalue of problem (\ref{ep0}).
Applying Lemma 2.1 to problem (\ref{ap}), there exists a
continuum $\mathscr{C}$ of nontrivial solutions of problem (\ref{ap}) emanating from $\left(\lambda_{1,\alpha},0\right)$.
Since (0,0) is the only solution of problem (\ref{ap}) for $\lambda = 0$, so
$\mathscr{C}\cap\left(\{0\}\times E\right)=\emptyset$.
By an argument similar to that of [\ref{DPM}, Lemma 3.1] with obvious changes, we may obtain that
\begin{equation}
{\mathscr{C}}\subseteq\left(\left\{\left(\lambda_{1,\alpha},0\right)\right\}\cup\left(\mathbb{R}\times \mathbb{S}^+\right)\cup\left(\mathbb{R}\times \mathbb{S}^-\right)\right)\nonumber
\end{equation}
and $\mathscr{C}$ is unbounded in $\mathbb{R}\times E$.

By Lemma 2.4, we know that $\mathscr{C}$ can be split into $\mathscr{C}^+$ and $\mathscr{C}^-$ such that $\mathscr{C}^+\cap\mathscr{C}^-=\left\{\left(\lambda_{1,\alpha},0\right)\right\}$.
Now, let us show that $\mathscr{C}^+\subseteq\left(\left\{\left(\lambda_{1,\alpha},0\right)\right\}\cup\left(\mathbb{R}\times \mathbb{S}^+\right)\right)$.
According to Lemma 2.2, there exists a bounded open neighborhood $\mathbb{O}$ of $\left(\lambda_{1,\alpha},0\right)$ such that
\begin{equation}
\left(\mathscr{C}^+\cap\mathbb{O}\right)\subseteq\left(\left\{\left(\lambda_{1,\alpha},0\right)\right\}\cup\left(\mathbb{R}\times \mathbb{S}^+\right)\right)\,\, \text{or}\,\, \left(\mathscr{C}^+\cap\mathbb{O}\right)\subseteq\left(\left\{\left(\lambda_{1,\alpha},0\right)\right\}\cup\left(\mathbb{R}\times \mathbb{S}^{-}\right)\right).\nonumber
\end{equation}
Without loss of generality, we assume that
$\left(\mathscr{C}^+\cap\mathbb{O}\right)\subseteq\left(\left\{\left(\lambda_{1,\alpha},0\right)\right\}\cup\left(\mathbb{R}\times \mathbb{S}^+\right)\right)$.
Suppose now the assertion is not true. Then $\mathscr{C}^+$  leaves $\mathbb{R}\times \mathbb{S}^+$ at some
point $\left(\lambda_*,u_*\right)\neq \left(\lambda_{1,\alpha},0\right)$. Necessarily $u_*\neq 0$, for otherwise $\lambda_*$ would be an eigenvalue of problem (\ref{ep0}) different to $\lambda_{1,\alpha}$ and an argument similar to the one of [\ref{DPM}, Lemma 3.1] would lead to a contradiction. A continuity argument shows that
$\left(\lambda_*,u_*\right)$ weakly satisfies
\begin{equation}
\left\{
\begin{array}{ll}
-\text{div}\left(\left\vert\nabla u_*\right\vert^{p-2}\nabla u_*\right)=\left(\lambda_* \alpha(x)+\lambda_*\frac{\xi\left(x,u_*\right)}{\varphi_p\left(u_*\right)}\right)\varphi_p\left(u_*\right)\,\,&\text{in}\,\,\Omega,\\
u_*=0 &\text{on}\,\,\partial\Omega.
\end{array}
\right.\nonumber
\end{equation}
Thus, $u_*\in C^{1,\delta}\left(\overline{\Omega}\right)$ and does not change sign in $\Omega$, say $u_*\geq 0$ in $\Omega$. By the strong maximum principle of [\ref{PS1}], we have $u_*>0$ in $\Omega$.
Obviously, $u_*$ cannot be approximated by any sequence of $\mathbb{S}^{-}$. This contradicts our assumption. Arguing in the same way leads to
$\mathscr{C}^-\subseteq\left(\left\{\left(\lambda_{1,\alpha},0\right)\right\}\cup\left(\mathbb{R}\times \mathbb{S}^-\right)\right)$.

Let us show that both $\mathscr{C}^+$ and $\mathscr{C}^-$ are unbounded. Consider
the following auxiliary problem
\begin{equation}
\left\{
\begin{array}{ll}
-\text{div}\left(\vert\nabla u\vert^{p-2}\nabla u\right)=\lambda \alpha(x)\varphi_p(u)+\lambda\widehat{\xi}(x,u)\,\,&\text{in}\,\,\Omega,\\
u=0 &\text{on}\,\,\partial\Omega,
\end{array}
\right.\nonumber
\end{equation}
where $\widehat{\xi}$ is defined by
\begin{equation}
\widehat{\xi}(x,s)=\left\{
\begin{array}{ll}
\xi(x,s)\,\,~~~~~&\text{if}\,\,s\geq 0,\\
-\xi(x,-s)\,\,&\text{otherwise}.
\end{array}
\right.\nonumber
\end{equation}
The previous
argument shows that an unbounded continuum $\widetilde{\mathscr{C}}$
bifurcates from $\left(\lambda_{1,\alpha},0\right)$ and can be split into $\widetilde{\mathscr{C}}^+$ and $\widetilde{\mathscr{C}}^-$
with $\widetilde{\mathscr{C}}^\nu$ connected, $\widetilde{\mathscr{C}}^\nu\subseteq\left(\left\{\left(\lambda_{1,\alpha},0\right)\right\}\cup\left(\mathbb{R}\times \mathbb{S}^\nu\right)\right)$.
It is easy to see that $\widetilde{\mathscr{C}}^-=-\widetilde{\mathscr{C}}^+$. It follows that both $\widetilde{\mathscr{C}}^+$ and $\widetilde{\mathscr{C}}^-$ are unbounded.
It is clear that $\widetilde{\mathscr{C}}^+\subseteq {\mathscr{C}}^+$. Therefore, ${\mathscr{C}}^+$ must be unbounded. A symmetric argument
shows that ${\mathscr{C}}^-$ is also unbounded.\qed\\

It is easy to verify that $\left(\lambda_{1,\alpha},0\right)$ is the unique bifurcation point of one-sign solutions of problem (\ref{ap}) from the trivial solution axis.
Now, we can prove Theorem 1.3 by virtue of Proposition 3.1.
\\ \\
\textbf{Proof of Theorem 1.3.} We claim that $u$ is also a solution of problem (\ref{zeros}) for any $(\lambda,u)\in \mathscr{C}$ which is obtained in Proposition 3.1. It is enough to show that $u\in[b(x),a(x)]$ for any $(\lambda,u)\in \mathscr{C}$. In fact if there exists $x_0\in\Omega$ such that $u\left(x_0\right)>a\left(x_0\right)$ or $u\left(x_0\right)<b\left(x_0\right)$, the definition of $\widetilde{f}$ will imply that $u\equiv 0$.
This is impossible.
So one has $u\in[b(x),a(x)]$ for any $(\lambda,u)\in \mathscr{C}$. It follows that $u$ is also a solution of problem (\ref{zeros})
for any $(\lambda,u)\in \mathscr{C}$.

Next we show that the projection of $\mathscr{C}^+$ on $\mathbb{R}$ is unbounded. It is sufficient to show that the set $\left\{(\lambda, u)\in  \mathscr{C}^+:\lambda\in (0, d]\right\}$ is bounded for any fixed
$d\in (0,+\infty)$. Suppose on the contrary that there exists $\left(\lambda_n,u_n\right)\in \mathscr{C}^+$, $n\in \mathbb{N}$, such that
$\lambda_n\rightarrow \mu\leq d$, $u_n\rightarrow +\infty$ as $n\rightarrow +\infty$. Let $w_n =u_n/\left\Vert u_n\right\Vert$.
Then we have that
\begin{equation}
w_n=R_p\left(\lambda_n \frac{{f}\left(x,u_n(x)\right)}{\left\Vert u_n\right\Vert^{p-1}}\right).\nonumber
\end{equation}
Clearly, we have that
\begin{equation}
{f}\left(x,u_n\right)\leq \max_{\overline{\Omega}\times [0,M]}\left\vert {f}(x,s)\right\vert,\nonumber
\end{equation}
where $M=\max_{x\in \overline{\Omega}} a(x)$.
It means that
\begin{equation}
\lambda_n \frac{f\left(x,u_n\right)}{\left\Vert u_n\right\Vert^{p-1}}\rightarrow 0\nonumber
\end{equation}
as $n\rightarrow +\infty$.
By the compactness of $R_p$, we obtain that for some convenient subsequence
$w_n\rightarrow w_0$ as $n\rightarrow+\infty$. Letting $n\rightarrow+\infty$,
we obtain that
$w_0\equiv0$. This contradicts $\left\Vert w_0\right\Vert=1$.
This together with the fact that $\mathscr{C}^+$ joins $\left(\lambda_{1,\alpha}, 0\right)$ to infinity yields that
\begin{equation}
\left(\lambda_{1,\alpha},+\infty\right)\subseteq \text{Proj}\left(\mathscr{C}^+\right),\nonumber
\end{equation}
where $\text{Proj}\left(\mathscr{C}^+\right)$ denotes the projection of $\mathscr{C}^+$ on $\mathbb{R}$.
Similarly, we can show that
\begin{equation}
\left(\lambda_{1,\alpha},+\infty\right)\subseteq \text{Proj}\left(\mathscr{C}^-\right).\nonumber
\end{equation}
\indent Now the existence of $u_0^+$ and $u_0^-$ can be got from Proposition 3.1 and the above argument immediately.
Moreover, we have $u_0^+\not\equiv a$ and $u_0^-\not\equiv b$ on $\overline{\Omega}$, that is to say, there exist some $x_0$, $x_1\in\overline{\Omega}$ such that $u_0^+\left(x_0\right)<a\left(x_0\right)$ and $u_0^-\left(x_1\right)>b\left(x_1\right)$. Otherwise, it contradicts hypothesis (H2).\qed\\

To prove Theorem 1.4, we propose the definition of linearly stable solution. For any $\phi\in E$ and nontrivial solution $u$ of problem (\ref{zeros}), the linearized problem of problem (\ref{zeros}) about $u$ at the direction $\phi$ is
\begin{equation}\label{lp}
\left\{
\begin{array}{ll}
-(p-1)\text{div}\left(\vert \nabla u\vert^{p-2}\nabla\phi\right)-\lambda f_u(x,u)\phi=\mu \phi\,\, &\text{in}
\,\,\Omega,\\
\phi=0&\text{on}\,\, \partial\Omega,
\end{array}
\right.
\end{equation}
where $f_u(x,u)\phi$ denotes the Fr\'{e}ch{e}t derivative of $f$ about $u$ at the direction $\phi$.
A solution $u$ of problem (\ref{zeros}) is stable if all eigenvalues of problem (\ref{lp})
are positive, otherwise it is unstable.
We define the Morse index $M(u)$ of a solution $u$ to problem (\ref{zeros})
to be the number of negative eigenvalues of problem (\ref{lp}).
A solution $u$ of problem (\ref{zeros}) is degenerate if $0$ is an eigenvalue of
problem (\ref{lp}), otherwise it is non-degenerate.
The following lemma is our main stability result for the one-sign solution.
\\ \\
\textbf{Lemma 3.1.} \emph{Let (H6) hold. Then any positive or negative solution $u$ of problem (\ref{zeros}) is stable,
hence, non-degenerate and Morse index $M(u)=0$.}
\\ \\
\textbf{Proof.} Without loss of generality, let $u$ be a positive solution of problem (\ref{zeros}), and let $\left(\mu_1, \phi_1\right)$ be
the corresponding principal eigen-pair of problem (\ref{lp}) with $\phi_1>0$ in $\Omega$. We
notice that $u$ and $\phi_1$ satisfy the problems
\begin{equation}\label{lp1}
\left\{
\begin{array}{ll}
-\text{div}\left(\varphi_p(\nabla u)\right)-\lambda f(x,u)=0\,\, &\text{in\,\,}\Omega,\\
u=0&\text{on}\,\, \partial\Omega
\end{array}
\right.
\end{equation}
and
\begin{equation}\label{lp2}
\left\{
\begin{array}{ll}
-(p-1)\text{div}\left(\vert \nabla u\vert^{p-2}\nabla\phi\right)-\lambda f_u(x,u)\phi=\mu \phi\,\, &\text{in}
\,\,\Omega,\\
\phi=0&\text{on}\,\, \partial\Omega.
\end{array}
\right.
\end{equation}
Multiplying the first equation of problem (\ref{lp2}) by $u$ and the first equation of problem (\ref{lp1}) by $(p-1)\phi_1$, subtracting and integrating, we obtain
\begin{equation}
\mu_1\int_\Omega \phi_1 u\,dx=\lambda\int_\Omega \phi_1\left((p-1)f(x,u)-f_u(x,u)u\right)\,dx.\nonumber
\end{equation}
By some simple computations, we can show that (H6) follows that $(p-1)f(x,s)-f_s(x,s)s\geq 0$ for any $s\geq0$.
Since $u> 0$ and $\phi_1> 0$ in $\Omega$, then $\mu_1 > 0$ and the positive solution
$u$ must be stable. \qed
\\ \\
\textbf{Proof of Theorem 1.4.} Define $F : \mathbb{R} \times E \rightarrow \mathbb{R}$ by
\begin{equation}
F(\lambda,u)=-\text{div}\left(\varphi_p(\nabla u)\right)-\lambda f(x,u).\nonumber
\end{equation}
From Lemma 3.1, we know that any one-sign solution $(\lambda, u)$ of problem (\ref{zeros}) is stable.
Therefore, at any one-sign solution $\left(\lambda^*, u^*\right)$, we can apply Implicit Function Theorem (see for example [\ref{De}]) to
$F(\lambda, u) = 0$, and all the solutions of $F(\lambda, u) = 0$ near $\left(\lambda^*, u^*\right)$ are on a curve
$(\lambda, u(\lambda))$ with $\left\vert \lambda-\lambda^*\right\vert\leq\varepsilon$ for some small $\varepsilon > 0$.
Furthermore, from Lemma 2.2, we can see that $\mathscr{C}=\mathscr{C}^+\cup \mathscr{C}^-$ near
$\left(\lambda_{1,\alpha},0\right)$ is given by a curve $(\lambda(s),u(s))=\left(\lambda_{1,\alpha}+o(1),s\varphi_{1,\alpha}+o(s)\right)$ for $s$ near $0$. Moreover, we can distinguish between two portions of this curve by $s\geq 0$ and $s< 0$. Therefore, the unbounded continua $\mathscr{C}^+$ and $\mathscr{C}^-$ are all curves.

To complete the proof, it suffices to show that $u^+(\lambda,\cdot)$ ($u^-(\lambda,\cdot)$) is increasing (decreasing) with respect to $\lambda$.
We only prove the case of $u^+(\lambda,\cdot)$. The proof of
$u^-(\lambda,\cdot)$ can be given similarly. Since $u^+(\lambda,\cdot)$ is differentiable with respect to
$\lambda$ (as a consequence of Implicit Function Theorem), then
$\frac{du^+(\lambda,\cdot)}{d\lambda}$ satisfies
\begin{equation}
-(p-1)\text{div}\left(\left\vert\nabla u^+\right\vert^{p-2}\nabla\frac{du^+}{d\lambda}\right)
=\lambda  f_u\left(x,u^+\right)\frac{d u^+}{d \lambda}
+f\left(x,u^+\right).\nonumber
\end{equation}
By an argument similar to that of Lemma 3.1, we can show that
\begin{equation}
\int_\Omega\left(\lambda \left(f_u\left(x,u^+\right)u^+-(p-1)f\left(x,u^+\right)\right)\frac{d u^+}{d \lambda}+f\left(x,u^+\right)u^+\right)\,dx=0.\nonumber
\end{equation}
(H2) implies $f(x,s)s\geq 0$ for any $s\in \mathbb{R}$. So we get $\left(f_u\left(x,u^+\right)u^+-(p-1)f\left(x,u^+\right)\right)\frac{d u^+}{d \lambda}\leq 0$.
While, (H6) shows that $f_u\left(x,u^+\right)u^+-(p-1)f\left(x,u^+\right)\leq 0$. Therefore, we have $\frac{d u^+}{d \lambda}\geq0$.\qed
\\
\\
\textbf{Remark 3.1.} From the proofs of Proposition 3.1, Theorem 1.3 and 1.4, we can see that all of these results can be obtained under less restrictive hypotheses than (H2). In fact, we only need that hypotheses (H2) holds on interval $[b,a]$. Moreover, signum condition is not needed in Proposition 3.1 and Theorem 1.3.

\section{Linear growth at infinity}
\bigskip
\quad\, In this section, we consider the case of linear growth of $f$ at infinity.
We first use Theorem 1.1 to study the bifurcation phenomenon of problem (\ref{zeros}) from infinity. Then we present the proofs of Theorem 1.5 and 1.6.
\\ \\
\textbf{Proposition 4.1.} \emph{Under hypotheses (H1), (H2) and (H4), the pair $\left(\lambda_{1,\beta},\infty\right)$ is a bifurcation point of problem (\ref{zeros}). Moreover, there exist two
unbounded continua of the set of nontrivial solutions of
problem (\ref{zeros}) in $\mathbb{R}\times E$ bifurcating
from $\left(\lambda_{1,\beta},\infty\right)$, $\mathscr{D}^+$ and $\mathscr{D}^-$, such that $\mathscr{D}^\nu\subseteq\left(\left\{\left(\lambda_{1,\beta},\infty\right)\right\}\cup\left(\mathbb{R}\times \mathbb{S}^\nu\right)\right)$.}
\\
\\
\textbf{Proof.} Let $\eta:\Omega\times\mathbb{R}\rightarrow \mathbb{R}$ be such that
\begin{equation}
f(x,s)=\beta(x)\varphi_p(s)+\eta(x,s)\nonumber
\end{equation}
with
\begin{equation}
\lim_{\vert s\vert\rightarrow+\infty}\frac{\eta(x,s)}{\varphi_p(s)}=0.\nonumber
\end{equation}
Let us consider
\begin{equation}
\left\{
\begin{array}{ll}
-\Delta_p u=\lambda \beta(x)\varphi_p(u)+\lambda\eta(x,u)\,\, &\text{in}\,\, \Omega,\\
u=0  &\text{on}\,\,\partial\Omega
\end{array}
\right.\nonumber
\end{equation}
as a bifurcation problem from infinity.

Obviously, problem (\ref{zeros}) can be equivalently written as
\begin{equation}\label{oei}
u=R_p\left(\lambda \widetilde{K}_pu+\widetilde{H}(\lambda,u)\right),
\end{equation}
where $\widetilde{K}_p(u)=\beta\varphi_p(u)$, $\widetilde{H}(\lambda,\cdot)$ denotes the usual Nemitsky operator associated with $\lambda\eta$.
Let
$L(\lambda)u=R_p\left(\lambda \widetilde{K}_pu\right)$ and ${H}(\lambda,u)=R_p\left(\lambda \widetilde{K}_pu+\widetilde{H}(\lambda,u)\right)-L(\lambda)u$.
Then it is easy to see that $L(\cdot):E\rightarrow E$ is homogeneous completely continuous.
From hypotheses (H1) and (H4),
the right hand side of equation (\ref{oei}) defines a completely continuous operator from $\mathbb{R}\times E$ into $E$.
So ${H}:\mathbb{R}\times E\rightarrow E$ is completely continuous.
Let
\begin{equation}
\widetilde{\eta}(x,u)=\max_{0\leq \vert s\vert\leq u}\vert \eta(x,s)\vert\,\,
\text{for any}\,\, x\in\Omega,\nonumber
\end{equation}
then $\widetilde{\eta}$ is nondecreasing with respect to $u$. Define
\begin{equation}
\overline{\eta}(x,u)=\max_{u/2\leq \vert s\vert\leq u}\vert \eta(x,s)\vert\,\,
\text{for any}\,\, x\in\Omega.\nonumber
\end{equation}
Then we can see that
\begin{equation}
\lim_{ u\rightarrow +\infty}\frac{\overline{\eta}(x,u)}{u^{p-1}}=0\,\,\text{and}\,\,\widetilde{\eta}(x,u)\leq \widetilde{\eta}\left(x,\frac{u}{2}\right)+\overline{\eta}(x,u).\nonumber
\end{equation}
It follows that
\begin{equation}
\limsup_{ u\rightarrow +\infty}\frac{\widetilde{\eta}(x,u)}{u^{p-1}}\leq \limsup_{ u\rightarrow +\infty}\frac{\widetilde{\eta}\left(x,\frac{u}{2}\right)}{u^{p-1}}
=\limsup_{ u/2\rightarrow +\infty}\frac{\widetilde{\eta}\left(x,\frac{u}{2}\right)}{2^{p-1}\left(\frac{u}{2}\right)^{p-1}}.\nonumber
\end{equation}
So we have
\begin{equation}\label{eg0+}
\lim_{ u\rightarrow +\infty}\frac{\widetilde{\eta}(x,u)}{u^{p-1}}=0.
\end{equation}
Further it follows from (\ref{eg0+}) that
\begin{equation}
\frac{\eta\left(x,u\right)}{\left\Vert u\right\Vert^{p-1}} \leq\frac{
\widetilde{\eta}\left(x,\left\vert u\right\vert\right)}{\left\Vert u\right\Vert^{p-1}} \leq \frac{
\widetilde{\eta}\left(x,\left\Vert u\right\Vert_\infty\right)}{\left\Vert u\right\Vert^{p-1}}\leq \frac{
\widetilde{\eta}\left(x,\left\Vert u\right\Vert\right)}{\left\Vert u\right\Vert^{p-1}}\rightarrow0\,\,  \text{as}\,\, \Vert u\Vert\rightarrow +\infty\nonumber
\end{equation}
uniformly in $x\in\Omega$.
It follows that ${H}=o(\Vert u\Vert)$ near $u=\infty$ uniformly on bounded $\lambda$ intervals.
By an argument similar to that of [\ref{DPM}, Proposition 2.2], we can get that
\begin{equation}
\deg_E\left(I-L(\lambda), B_r(0),0 \right)=\left\{
\begin{array}{ll}
1 \,\, &\text{if}\,\, 0<\lambda<\lambda_{1,\beta},\\
-1\,\,&\text{if}\,\,\lambda_{1,\beta}<\lambda<\lambda_{2,\beta},
\end{array}
\right.\nonumber
\end{equation}
where $\lambda_{2,\beta}$ denotes the second eigenvalue of problem (\ref{ep0}) with the weighted function $\beta$.

Applying Proposition 2.1 to problem (\ref{oei}), there exists a
continuum $\mathscr{D}$ of solutions of problem (\ref{zeros}) meeting $\left(\lambda_{1,\beta},\infty\right)$ and satisfying at least one of the alternatives of Proposition 2.1.
We claim that 1$^o$ of Proposition 2.1 does not occur. Otherwise, $\mathscr{D}$ will crosse the line $\mathbb{R}\times \{a\}$ or $\mathbb{R}\times \{b\}$ in $\mathbb{R}\times E$.
That is to see $a$ or $b$ is a solution of problem (\ref{zeros}) which contradicts the assumption of (H2).
So 2$^o$ of Proposition 2.1 occurs.

We claim that $\mathscr{D}-\mathscr{M}$ has an unbounded projection on $\mathbb{R}$.
By Theorem 1.1, we know that $\mathscr{D}$ can be split into $\mathscr{D}^+$ and $\mathscr{D}^-$ such that $\mathscr{D}^+\cap\mathscr{D}^-=\left\{\left(\lambda_{1,\beta},\infty\right)\right\}$.
According to Theorem 1.1, there exists a bounded open neighborhood $\mathscr{N}$ of $\left(\lambda_{1,\beta},\infty\right)$ such that
$\left(\mathscr{D}^\nu\cap\mathscr{N}\right)\subseteq\left(\left\{\left(\lambda_{1,\beta},\infty\right)\right\}\cup\left(\mathbb{R}\times \mathbb{S}^\nu\right)\right)$ or $\left(\mathscr{D}^\nu\cap\mathscr{N}\right)\subseteq\left(\left\{\left(\lambda_{1,\beta},\infty\right)\right\}\cup\left(\mathbb{R}\times \mathbb{S}^{-\nu}\right)\right)$. Without loss of generality, we assume that
$\left(\mathscr{D}^\nu\cap\mathscr{N}\right)\subseteq\left(\left\{\left(\lambda_{1,\beta},\infty\right)\right\}\cup\left(\mathbb{R}\times \mathbb{S}^\nu\right)\right)$. Suppose on the contrary that $\left(\mathscr{D}^+-\left(\mathscr{D}^+\cap\mathscr{N}\right)\right)\not\subset \mathbb{R}\times\mathbb{S}^+$, then $\mathscr{D}^+$ leaves
$\mathbb{R}\times \mathbb{S}^+$ at some point $(\mu, u)\neq \left(\lambda_{1,\beta},\infty\right)$ . If $u\equiv 0$, i.e., 1$^o$ of Proposition 2.1 occurs,  which is a contradiction. So $u\not \equiv0$. A continuity argument shows that
$\left(\mu,u\right)$ weakly satisfies
\begin{equation}
\left\{
\begin{array}{ll}
-\text{div}\left(\left\vert\nabla u\right\vert^{p-2}\nabla u\right)=\mu f\left(x,u\right)\,\,&\text{in}\,\,\Omega,\\
u=0&\text{on}\,\,\partial\Omega.
\end{array}
\right.\nonumber
\end{equation}
Thus, $u\in C^{1,\delta}\left(\overline{\Omega}\right)$ and does not change sign in $\Omega$, say $u\geq 0$ in $\Omega$. So $u>0$ in $\Omega$ by the strong maximum principle of [\ref{PS1}].
This contradicts our assumption. Thus, we have $\mathscr{D}^+\subseteq\left(\left\{\left(\lambda_{1,\beta},\infty\right)\right\}\cup\left(\mathbb{R}\times \mathbb{S}^+\right)\right)$. Arguing in the same way yields to
$\mathscr{D}^-\subseteq\left(\left\{\left(\lambda_{1,\beta},\infty\right)\right\}\cup\left(\mathbb{R}\times \mathbb{S}^-\right)\right)$.
Therefore, we have that
\begin{equation}\label{oei01}
{\mathscr{D}^\nu}\subseteq\left(\left\{\left(\lambda_{1,\beta},\infty\right)\right\}\cup\left(\mathbb{R}\times \mathbb{S}^\nu\right)\right).
\end{equation}
Now we show that the case of $\mathscr{D}-\mathscr{M}$ meeting $\lambda_{j,\beta}\times \{\infty\}$ for some $j>1$ is impossible, where $\lambda_{j,\beta}$ denotes the $j$th eigenvalue of problem (\ref{ep0}) with the weighted function $\beta$.
Assume on the contrary that $\mathscr{D}-\mathscr{M}$ meets $\lambda_{j,\beta}\times \{\infty\}$ for some $j>1$. So there exists a neighborhood $\widetilde{\mathscr{N}}\subset\widetilde{\mathscr{M}}$ of $\lambda_{j,\beta}\times \{\infty\}$ such that $u$ must change sign for any $(\lambda,u)\in\left(\mathscr{D}-\mathscr{M}\right)\cap\left(\widetilde{\mathscr{N}}\setminus\left(\lambda_{j,\beta}\times \{\infty\}\right)\right)$, where $\widetilde{\mathscr{M}}$ is a neighborhood of $\lambda_{j,\beta}\times \{\infty\}$ which satisfies the assumptions of Proposition 2.1. This contradicts (\ref{oei01}). Thus, Proposition 2.1 shows the desired conclusion.

Finally, we will show that both $\mathscr{D}^+$ and $\mathscr{D}^-$ have unbounded projections on $\mathbb{R}$. Consider
the following auxiliary problem
\begin{equation}
\left\{
\begin{array}{ll}
-\text{div}\left(\vert\nabla u\vert^{p-2}\nabla u\right)=\lambda \beta(x)\varphi_p(u)+\lambda\eta{\xi}(x,u)\,\,&\text{in}\,\,\Omega,\\
u=0 &\text{on}\,\,\partial\Omega,
\end{array}
\right.\nonumber
\end{equation}
where $\widetilde{\eta}$ is defined by
\begin{equation}
\widehat{\eta}(x,s)=\left\{
\begin{array}{ll}
\eta(x,s)\,\,~~~~~&\text{if}\,\,s\geq 0,\\
-\eta(x,-s)\,\,&\text{otherwise}.
\end{array}
\right.\nonumber
\end{equation}
The previous
argument shows that an unbounded continuum $\widetilde{\mathscr{D}}$
bifurcates from $\left(\lambda_{1,\beta},\infty\right)$ and can be split into $\widetilde{\mathscr{D}}^+$ and $\widetilde{\mathscr{D}}^-$
with $\widetilde{\mathscr{D}}^\nu$ connected, $\widetilde{\mathscr{D}}^\nu\subseteq\left(\left\{\left(\lambda_{1,\beta},\infty\right)\right\}\cup\left(\mathbb{R}\times \mathbb{S}^\nu\right)\right)$.
It is easy to see that $\widetilde{\mathscr{D}}^-=-\widetilde{\mathscr{D}}^+$. It follows that both $\widetilde{\mathscr{D}}^+$ and $\widetilde{\mathscr{D}}^-$ have unbounded projections on $\mathbb{R}$.
It is clear that $\widetilde{\mathscr{D}}^+\subseteq {\mathscr{D}}^+$. Therefore, ${\mathscr{D}}^+$ has an unbounded projection on $\mathbb{R}$. A symmetric argument
shows that ${\mathscr{D}}^-$ also has an unbounded projection on $\mathbb{R}$.
\qed\\

It is easy to verify that $\left(\lambda_{1,\beta},\infty\right)$ is the unique bifurcation point of one-sign solutions of problem (\ref{zeros}) from infinity. Now, we can prove Theorem 1.5 and 1.6 by using of Proposition 4.1.
\\ \\
\textbf{Proof of Theorem 1.5.} The conclusion (i) is a direct corollary of Proposition 4.1. So we only need to show the conclusion (ii).
Now, noting Proposition 3.1 and 4.1, it is sufficient to show that $\mathscr{C}^\nu\cap \mathscr{D}^\nu=\emptyset$. Otherwise, $\mathscr{D}^\nu$ meets $\mathscr{R}=\{(\lambda,0):\lambda\in \mathbb{R}\}$. Hence, it will crosse the line $\mathbb{R}\times \{a\}$ or $\mathbb{R}\times \{b\}$ in $\mathbb{R}\times E$.
That is to see $a$ or $b$ is a solution of problem (\ref{zeros}) which contradicts hypothesis (H2).\qed\\
\\
\textbf{Proof of Theorem 1.6.} The proof is similar to that of Theorem 1.5. So we omit it here.\qed

\section{Superlinear growth at infinity}
\bigskip
\quad\, In this section, we consider the case of superlinear growth of $f$ at infinity.
We first use Theorem 1.2 and Proposition 4.1 to study the bifurcation phenomenon of problem (\ref{zeros}) from infinity. Then we give the proof of Theorem 1.7. Finally, we study the asymptotic behavior of the solutions when $\lambda\rightarrow+\infty$ and prove Theorem 1.8.
\\ \\
\textbf{Proposition 5.1.} \emph{Under hypotheses (H1)--(H3) and (H5), the pair $\left(0,\infty\right)$ is a bifurcation point of problem (\ref{zeros}). Moreover, there exist two
unbounded continua of the set of nontrivial solutions of
problem (\ref{zeros}) in $\mathbb{R}\times E$ bifurcating
from $\left(0,\infty\right)$, $\mathscr{D}^+$ and $\mathscr{D}^-$, such that $\mathscr{D}^\nu\subseteq\left(\left\{\left(0,\infty\right)\right\}\cup\left(\mathbb{R}\times \mathbb{S}^\nu\right)\right)$.}
\\
\\
\textbf{Proof.} For any $n\in \mathbb{N}$, define
\begin{eqnarray}
f_n(x,s)=\left\{
\begin{array}{lll}
f(x,s),\,\, & s\in\left[-n,n\right],\\
\frac{\beta(x)\varphi_p(2n)-f(x,n)}{n}(s-n)+f(x,n), & s\in \left(\left(-2n,-n\right)\cup
\left(n,2n\right)\right),\\
n\beta(x)\varphi_p(s), & s\in\left((-\infty,-2n]\cup\left[2n,+\infty\right)\right).
\end{array}
\right.\nonumber
\end{eqnarray}
Now, consider the following problem
\begin{equation}\label{4.2}
\left\{
\begin{array}{ll}
-\Delta_p u=\lambda f_n(x,u)\,\, &\text{in}\,\, \Omega,\\
u=0  &\text{on}\,\,\partial\Omega.
\end{array}
\right.
\end{equation}
Clearly, we can see that $\lim_{n\rightarrow+\infty}f_n(x,s)=f(x,s)$ and
\begin{equation}
\lim_{\vert s\vert\rightarrow +\infty}\frac{f_n(x,s)}{\varphi_p(s)}=n\beta(x)\,\,\text{uniformly in}\,\, x\in \Omega.\nonumber
\end{equation}
Proposition 4.1 implies
that there exist two sequence unbounded continua $\mathscr{D}_n^\nu$ of solution set of problem (\ref{4.2}) emanating from $\left(\lambda_{1,\beta}/n, \infty\right)$.

Taking $z^*=\left(0, \infty\right)$, clearly $z^*\in \liminf_{n\rightarrow +\infty}\mathscr{D}_n^\nu$ with $\left\Vert z^*\right\Vert_{\mathbb{R}\times E}=+\infty$.
Define mapping $T:X\rightarrow X$ such that
\begin{equation}
T(\lambda,u)=\left\{
\begin{array}{ll}
\left(\lambda,\frac{u}{\Vert u\Vert^2}\right)\,\, &\text{if}\,\, 0<\Vert u\Vert<+\infty,\\
(\lambda,0)  &\text{if}\,\,\Vert u\Vert=+\infty,\\
(\lambda,\infty)  &\text{if}\,\,\Vert u\Vert=0.
\end{array}
\right.\nonumber
\end{equation}
It is easy to verify that $T$ is homeomorphism and $\left\Vert T\left( z^*\right)\right\Vert_{\mathbb{R}\times E}=0$. The compactness of $R_p$ implies that $\left(\cup_{n=1}^{+\infty} T\left(\mathscr{D}_n^\nu\right)\right)\cap B_R$ is pre-compact. Theorem 1.2 implies
that $\mathscr{D}^\nu=\limsup_{n\rightarrow +\infty}\mathscr{D}_n^\nu$ is unbounded closed
connected such that $z^*\in \mathscr{D}^\nu$.

Next we show that the projection of $\mathscr{D}^\nu$ on $\mathbb{R}$ is unbounded. From the argument of Proposition 4.1, we have known that
$\mathscr{D}_n^\nu$ has unbounded projections on $\mathbb{R}$. By Proposition 2.3, for each $\epsilon> 0$ there exists an $m$ such that for every $n > m$, $\mathscr{D}_n^\nu\subset V_\epsilon\left(\mathscr{D}^\nu\right)$.
This implies that
\begin{equation}
\left(\frac{\lambda_{1,\beta}}{n},+\infty\right)\subseteq \text{Proj}\left(\mathscr{D}_n^\nu\right)\subseteq \text{Proj}\left(V_\epsilon\left(\mathscr{D}^\nu\right)\right).\nonumber
\end{equation}
It follows that the projection of $\mathscr{D}^\nu$ is unbounded on $\mathbb{R}$. Therefore, we have that
\begin{equation}
\left(0,+\infty\right)= \text{Proj}\left(\mathscr{D}^\nu\right).\nonumber
\end{equation}
\indent Finally, we show $\mathscr{D}^\nu\subseteq\left(\left\{\left(0,\infty\right)\right\}\cup\left(\mathbb{R}\times \mathbb{S}^\nu\right)\right)$.
Proposition 4.1 shows that
\begin{equation}
\mathscr{D}_n^\nu\subseteq\left(\left\{\left(\frac{\lambda_{1,\beta}}{n},\infty\right)\right\}\cup\left(\mathbb{R}\times \mathbb{S}^\nu\right)\right).\nonumber
\end{equation}
It follows that $\nu u$ is a nonnegative solution of problem (\ref{zeros}) for any $(\lambda,u)\in \left(\mathscr{D}^\nu\setminus \left\{\left(0,\infty\right)\right\}\right)$. The signum condition and the strong maximum principle [\ref{PS1}] imply that $\nu u>0$ in $\Omega$.
So we have $\mathscr{D}^\nu\subseteq\left(\left\{\left(0,\infty\right)\right\}\cup\left(\mathbb{R}\times \mathbb{S}^\nu\right)\right)$.
\qed\\
\\
\textbf{Proof of Theorem 1.7.} In view of Proposition 5.1, an argument similar to that of Theorem 1.5 can yield the desired conclusions. \qed
\\

In order to prove Theorem 1.8, we need the following priori estimate for the solution of problem (\ref{zeros}).
\\ \\
\textbf{Lemma 5.1.} \emph{Suppose assumptions (H1), (H3) and ($\widetilde{H}$5) hold. Given $\widetilde{\lambda}>0$, there exists a constant $D$ such that if $u\in E$ is an one-sign solution of problem (\ref{zeros}) with $\lambda>\widetilde{\lambda}$, then
\begin{equation}
\Vert u\Vert_\infty\leq D,\nonumber
\end{equation}
where $D$ only depends on $\widetilde{\lambda}$.}
\\ \\
\textbf{Proof.} The conclusion involving the case of $u$ being positive is just Lemma 4.1 of [\ref{IMSU}]. Moreover, the negative case is similar to the positive case. \qed\\

Furthermore, we also need the following Liouville-type theorem
for a non-positive function with zero.
\\ \\
\textbf{Lemma 5.2.} \emph{Let $h :(-\infty,0]\rightarrow(-\infty,0]$ be a continuous function satisfying the following four assumptions:}\\

(h1) \emph{There exists a constant $\overline{b} < 0$ such that $h(s)=0$ if $s=0$ or $s=\overline{b}$, $h(s)s>0$ for $s<0$, $\neq \overline{b}$.}

(h2) \emph{There exist constants $\gamma > 0$ and $q\in \left(p, p_*\right)$ such that $h(s) \leq-\gamma \left(\overline{b}-s\right)^{q-1}$ for $s<\overline{b}$}.

(h3) \emph{There exists a constant $\overline{\alpha} > 0$ such that $\lim_{s\rightarrow0^-} h(t)/\varphi_p(s)=\overline{\alpha}$.}

(h4) \emph{There exists a constant $\Lambda> 0$ such that $-\Lambda\left(\vert s\vert^{q-1}+1\right)\leq h(s)$ for $s\leq 0$.}
\\

\noindent \emph{Then any $C^1$ weak solution of the problem}
\begin{equation}
\left\{
\begin{array}{l}
-\Delta_p v=h(v) \,\, \text{in}\,\, \mathbb{R}^N,\\
v\leq 0,
\end{array}
\right.\nonumber
\end{equation}
\emph{is either the constant function $w\equiv 0$, or else $w\equiv \overline{b}$.}\\
\\
\textbf{Proof.} Let $w=-v$ and $f(t)=-h(-t)$. Applying Theorem 1.8 of [\ref{IMSU}] to
\begin{equation}
\left\{
\begin{array}{l}
-\Delta_p w=f(w) \,\, \text{in}\,\, \mathbb{R}^N,\\
w\geq 0,
\end{array}
\right.\nonumber
\end{equation}
the desired conclusion can be got immediately. \qed
\\ \\
\textbf{Proof of Theorem 1.8.} Noting Lemma 5.1 and 5.2, by an argument essentially same as that of [\ref{IMSU}, Theorem 1.7], we can get the desired conclusions. \qed


\end{document}